\documentclass[12pt,a4paper]{amsart2020}
\usepackage{amsfonts}
\usepackage[top=35mm, bottom=35mm, left=30mm, right=30mm]{geometry}
\usepackage[colorlinks=true,citecolor=blue]{hyperref}
\usepackage{mathptmx}
\usepackage{eucal}
\usepackage{graphicx}
\usepackage{mathrsfs}
\usepackage{amssymb}
\usepackage{amsmath}
\usepackage{amsthm}
\usepackage{tikz,pgfplots,float}

\usepackage{xcolor}
\newtheorem{thm}{Theorem}[section]

\newtheorem{lem}[thm]{Lemma}
\newtheorem{prop}[thm]{Proposition}

\theoremstyle{definition}
\newtheorem{defn}[thm]{Definition}

\newtheorem{rem}[thm]{Remark}
\numberwithin{equation}{section}

\newcommand{\N}{\mathbb{N}}
\newcommand{\F}{\mathcal{F}}

\newcommand{\calM}{\mathcal{M}}
\newcommand{\calR}{\mathcal{R}}
\newcommand{\calP}{\mathcal{P}}
\newcommand{\calF}{\mathcal{F}}
\newcommand{\calB}{\mathcal{B}}
\newcommand{\calE}{\mathcal{E}}
\newcommand{\calI}{\mathcal{I}}

\newcommand{\R}{\mathbb{R}}
\newcommand{\scrP}{\mathscr{P}}
\newcommand{\vp}{\varepsilon}
\DeclareMathOperator{\md}{d}

\begin{document}
\title[Variational principles for topological pressures on subsets]{Variational principles for topological pressures on subsets}
\author{Xing-fu Zhong}
\address{School of Mathematics and Statistics, Guangdong University of Foreign Studies\\
Guangzhou, 510006 P. R. China}
\email{xfzhong@gdufs.edu.cn}

\author{Zhi-jing Chen}
\address{School of Mathematics and Systems Science,
	Guangdong Polytechnic Normal University,
	Guangzhou 510665, P.R. China}
\email{chzhjing@mail2.sysu.edu.cn}


\begin{abstract}
In this paper, we investigate the relations between various types of topological pressures and different versions of measure-theoretical pressures. We extend Feng- Huang's variational principle for packing entropy to packing pressure and obtain two new variational principles for Pesin-Pitskel and packing pressures respectively. We show that various types of Katok pressures for an ergodic measure with respect to a potential function are equal to the sum of measure-theoretic entropy of this measure and the integral of the potential function. Moreover, we obtain Billingsley type theorem for packing pressure, which indicates that packing pressure can be determined by measure-theoretic upper local pressure of measures, and a variational principle for packing pressure of the set of generic points for any invariant ergodic Borel probability measure.
\end{abstract}


\maketitle

\section{Introduction}
Throughout this paper, a \emph{topological dynamical system} (TDS for short) is a pair $(X,T)$, where $X$ is a compact metric space with a metric $d$ and $T:X\to X$ is a continuous map. Let ${\calM}(X)$, $\calM(X,T)$, and $\calE(X,T)$ denote respectively the set of Borel probability measures, $T$-invariant Borel probability measures, and $T$-invariant ergodic Borel probability measures on $X$. By a \emph{Borel measure theoretical dynamical system}$(X,\calB(X),\mu,T)$ we mean $(X,\calB(X),\mu)$ is a Borel measure space and $T$ is a Borel measure preserving transformation.

Kolmogorov~\cite{Kolmogorov1958A} introduced measure-theoretical entropy $h_\mu(T)$ for any $(X,\calB(X),\mu,T)$. Later, Adler, Konheim, and McAndrew~\cite{Adler1965} introduced topological entropy $h_{top}(T)$ for any TDS $(X,T)$. Dinaburg~\cite{Dinaburg1970} and Bowen~\cite{Bowen1971} independently gave equivalent definitions for topological entropy by using separating and spanning sets. The variational principle (see~\cite[Thm. 8.6]{Walters1982}) reveals the basic relation between topological entropy and measure-theoretical entropy: if $(X,T)$ is a TDS then
\[h_{top}(T)=\sup\{h_\mu(T):\mu\in\calM(X,T)\}.\]

Bowen~\cite{Bowen1973Topological} presented a type of topological entropy $h_{top}^B(Z,T)$ for any set $Z$ in a TDS $(X,T)$ in a way resembling Hausdorff dimension, which is the so-called Bowen topological entropy. Particularly, he showed that $h_{top}^B(X,T)=h_{top}(T)$ for any TDS $(X,T)$.

Topological pressure, as a non-trivial extension of topological entropy, was first introduced by Ruelle~\cite{Ruelle1973statistical} and extended to compact spaces with continuous transformations by Walters~\cite{Walters1982}. There is also a variational principle of topological pressure (see Thm.9.10 in~\cite{Walters1982}): if $(X,T)$ is a TDS then $P(T,f)=\sup\{h_\mu(T)+\int f\md\mu|\mu\in\calM(X,T)\},$ where $f$ is a continuous functions of $X$, and $P(T,f)$ is the pressure of $T$ with respect to $f$.
Pesin and Pitskel~\cite{Pesin1984Topological} further extended Bowen's results~\cite{Bowen1973Topological} to topological pressure and introduced a type of topological pressure $P^B(Z,T,f)$ for any set $Z$ in a TDS $(X,T)$, which we call Pesin-Pitskel topological pressure. They in~\cite{Pesin1984Topological} showed that, among other things, if $(X,T)$ is a TDS then $P(T,f)=P^B(X,T,f)$.

Inspired by the variational relation between topological entropy and measure-theoretical entropy, Feng and Huang~\cite{feng2012variational} introduced measure-theoretical lower and upper entropies and packing topological entropy, and they obtained two variational principles for Bowen entropy and packing entropy: if $Z\subset X$ is nonempty and compact then
\begin{align*}
h_{top}^B(Z,T)&=\sup\{\underline{h}_\mu(T):\mu\in\calM(X),~\mu(Z)=1\},\\
h_{top}^P(Z,T)&=\sup\{\overline{h}_\mu(T):\mu\in\calM(X),~\mu(Z)=1\},
\end{align*}
where $h_{top}^P(Z,T)$, $\underline{h}_\mu(T)$, and $\overline{h}_\mu(T)$ denote respectively the packing topological entropy of $Z$, measure-theoretical lower and upper entropies of $\mu$. Since then, Feng-Huang's variational principles have been extended to different systems and topological pressures; we refer the reader to \cite{Wang2012Variational,Xu2018Variational,Dou2017Topological,
Zhong2020Variationalfree,Tang2015Variational,
Zheng2016Bowen,Kong2014Slow,Liang2021Packing,Huang2020A} for more details. Tang, Cheng, and Zhao~\cite{Tang2015Variational} generalized Feng-Huang's variational principle of Bowen topological entropy to Pesin-Pitskel topological pressure: if $Z\subset X$ is nonempty and compact then
\[P^B(Z,T,f)=\sup\{\underline{P}_\mu(T,f):\mu\in\calM(X),~\mu(Z)=1\},\]
where $P^B(Z,T,f)$ denotes the Pesin-Pitskel pressure of $Z$ (see Definition~\ref{defn:pressures}) and $\underline{P}_\mu(T,f)$ denotes the measure-theoretical lower pressure of $\mu$ (see Definition~\ref{defn:M-L-U-pressures}).
Recently, Wang in~\cite{Wang2021some} obtained two new variational principles for Bowen and packing topological entropies by introducing Bowen entropy and packing entropy of measures in the sense of Katok. He showed that if $Z\subset X$ is nonempty and compact then
\begin{align*}
h_{top}^B(Z,T)&=\sup\{P_{\mu}^{KB}(T,0):\mu\in\calM(X),~\mu(Z)=1\},\\
h_{top}^P(Z,T)&=\sup\{P_{\mu}^{KP}(T,0):\mu\in\calM(X),~\mu(Z)=1\},
\end{align*}
where $P_{\mu}^{KB}(T,0)$ and $P_{\mu}^{KP}(T,0)$ denote respectively the Bowen and packing topological entropies of $\mu$ in the sense of Katok (see Definition~\ref{defn:M-K-pressures}). Wang and Chen~\cite{Wang2012Variational} studied the variational principles of BS dimensions and introduced packing topological pressure. Their work gives us a motivation for studying variational principles of Pesin-Pitskel and packing topological pressures.

Let $C(X,\R)$ denote the set of all continuous functions of $X$, and let $P_\mu^{B}(T,f)$ and $P_\mu^{KB}(T,f)$ denote respectively the Pesin-Pitskel pressure of $\mu$ (see Definition~\ref{defn:M-pressures}) and the Pesin-Pitskel pressure of $\mu$ in the sense of Katok (see Definition~\ref{defn:M-K-pressures}).

\begin{thm}\label{thm:variational-p-for-BP}
Let $(X,T)$ be a TDS, $f\in C(X,\R)$ and $Z\subset X$ be a nonempty compact set. Then
\begin{align*}
P^B(Z,T,f)&=\sup\{P_\mu^{B}(T,f):\mu\in{\calM}(X),~\mu(Z)=1\}\\
&=\sup\{P_\mu^{KB}(T,f):\mu\in{\calM}(X),~\mu(Z)=1\}.
\end{align*}
\end{thm}

\begin{rem}
(1) In fact, $P_\mu^{B}(T,f)=P_\mu^{KB}(T,f)$; see Proposition~\ref{prop:MKPs-equal-MPs}.

(2) When $f=0$, (1) of Theorem~1.2 in~\cite{Wang2021some} is a special case of Theorem~\ref{thm:variational-p-for-BP}.
\end{rem}

Let $P^P(Z,T,f)$, $\overline{P}_{\mu}(T,f)$, $P_{\mu}^{P}(T,f)$, and $P_{\mu}^{KP}(T,f)$ denote respectively the packing topological pressure of $Z$ (see Definition~\ref{defn:pressures}), measure-theoretical upper pressure of $\mu$ (see Definition~\ref{defn:M-L-U-pressures}), packing pressure of $\mu$ (see Definition~\ref{defn:M-pressures}), and packing pressure of $\mu$ in the sense of Katok (see Definition~\ref{defn:M-K-pressures}).

\begin{thm}\label{thm:variational-p-for-P}
Let $(X,T)$ be a TDS, $f\in C(X,\R)$, and $Z\subset X$ be a nonempty compact set. If $P^P(Z,T,f)>\|f\|_{\infty}$, where $\|f\|_{\infty}:=\sup_{x\in X}f(x)$, then
\begin{align*}
P^P(Z,T,f)&=\sup\{\overline{P}_{\mu}(T,f):\mu\in{\calM}(X),~\mu(Z)=1\}\\
&=\sup\{P_{\mu}^{P}(T,f):\mu\in{\calM}(X),~\mu(Z)=1\}\\
&=\sup\{P_{\mu}^{KP}(T,f):\mu\in{\calM}(X),~\mu(Z)=1\}.
\end{align*}
\end{thm}

\begin{rem}
(1) We will see from Proposition~\ref{prop:MKPs-equal-MPs} that $P_\mu^{P}(T,f)=P_\mu^{KP}(T,f)$.

(2) We obtain (\romannumeral1) of Theorem 1.3 in~\cite{feng2012variational} and (2) of Theorem 1.2 in~\cite{Wang2021some} as special cases of Theorem~\ref{thm:variational-p-for-P} by letting $f=0$.
\end{rem}

From Theorem 1.1 in~\cite{Katok1980Lyapunov}, we see that if $(X,T)$ is a TDS, $T$ is a homeomorphism of $X$, and $\mu\in\calE(X,T)$ then $h_\mu(T)=\underline{CP}^K_{\mu}(T,0)=\overline{CP}^K_{\mu}(T,0)$, where $\underline{CP}^K_{\mu}(T,0)$ and $\overline{CP}^K_{\mu}(T,0)$ are respectively the lower and upper capacity topological entropies of $\mu$ in the sense of Katok (see Definition~\ref{defn:M-K-pressures}). Theorem 11.6 in~\cite{Pesin1997Dimension} tells us that under the same assumptions and $f\in C(X,\R)$ then $P_\mu^{B}(T,f)=\underline{CP}_{\mu}(T,f)=\overline{CP}_{\mu}(T,f)=h_{\mu}(T)+\int f\md\mu$, where $\underline{CP}_{\mu}(T,f)$ and $\overline{CP}_{\mu}(T,f)$ are respectively the lower and upper capacity topological pressures of $\mu$ (see Definition~\ref{defn:M-pressures}). Wang recently in~\cite{Wang2021some} showed that $P_\mu^{KB}(T,0)=\underline{CP}^K_{\mu}(T,0)=P_{\mu}^{KP}(T,0)=\overline{CP}^K_{\mu}(T,0)=h_\mu(T)$.
We in this paper shall extend Wang's result to pressures in the sense of Katok.

\begin{thm}\label{thm:vp-among-mp}
Let $(X,T)$ be a TDS and $f\in C(X,\R)$. If $T$ is a homeomorphism of $X$ and
$\mu$ is a non-atomic Borel ergodic measure on $X$, then
\begin{align*}
P_\mu^{KB}(T,f)&=\underline{CP}^K_{\mu}(T,f)=P_{\mu}^{KP}(T,f)\\
&=\overline{CP}^K_{\mu}(T,f)=P_\mu^{B}(T,f)=\underline{CP}_{\mu}(T,f)\\
&=P_{\mu}^{P}(T,f)=\overline{CP}_{\mu}(T,f)=h_{\mu}(T)+\int f\md\mu.
\end{align*}
\end{thm}

Ma and Wen showed that Bowen entropy can be determined via the local entropies of measures in~\cite[Theorem 1]{Ma-Wen2008A-Billingsley}, which is an analogue of Billingsley's Theorem for the Hausdorff dimension. Tang, Cheng, and Zhao in~\cite{Tang2015Variational} extended this result to Pesin-Pitskel topological pressure. We have the following Billingsley type theorem for packing pressure.

\begin{thm}\label{thm:Billingsley-theorem}
Let $(X,T)$ be a TDS, $f\in C(X,\R)$, $\mu\in{\calM}(X)$, and $Z\subset X$. For $s\in\R$, the following properties hold:
\begin{itemize}
  \item[(1)] if $\overline{P}_{\mu}(x,T,f)\leq s$ for all $x\in Z$, then $P^P(Z,T,f)\leq s$;
  \item[(2)] if $\overline{P}_{\mu}(x,T,f)\geq s$ for all $x\in Z$ and $\mu(Z)>0$, then $P^P(Z,T,f)\geq s$,
\end{itemize}
where $\overline{P}_{\mu}(x,T,f)$ denotes the measure-theoretic upper local pressure of $x$ (see Section~\ref{sec:mea-pre}).
\end{thm}

Let $\mu\in\calM(X,T)$, $n\in\N$, and $x\in X$, define
\[\Upsilon_n(x)=\frac{1}{n}\sum_{i=0}^{n-1}\delta_{T^i(x)}.\]
Recall that $x$ is said to be \emph{generic for $\mu$} if $\lim_{n\to\infty}\Upsilon_n(x)=\mu$ in the weak$^\ast$-topology. Let $G_\mu$ denote the set of all generic points for $\mu$.

Bowen in~\cite{Bowen1973Topological} presented a variational principle of Bowen entropy for generic sets: if $\mu\in\calE(X,T)$ then $h_{top}^B(G_\mu,T)=h_\mu(T)$. Pesin and Pitskel in~\cite{Pesin1984Topological} extended this result to topological pressure: if $\mu\in\calE(X,T)$ then $P^B(G_\mu,T,f)=h_\mu(T)+\int f\md\mu$.
Zheng and Chen~\cite{Zheng2018Topological} generalized the variational principle of Bowen entropy of generic sets to countable discrete infinite amenable groups. Dou, Zheng, and Zhou~\cite{Dou2020Packing} obtained a variational principle of packing topological entropy of generic sets for amenable group actions. Wang et al.~\cite{Wang2020Bowen} gave a variational principle of Bowen entropy of generic sets for free flows without fixed points. Recently, Wang~\cite{Wang2021some} obtained a variational principle of packing topological entropy of generic sets for classical dynamical systems. We shall generalize Wang's result to packing topological pressure.

\begin{thm}\label{thm:va-pp-gene}
Let $(X,T)$ be a TDS, $\mu\in\calE(X,T)$, and $f\in C(X,\R)$. Then
\[P^P(G_\mu,T,f)=h_u(T)+\int f\md\mu.\]
\end{thm}

\begin{rem}
When $f=0$, we obtain Theorem 1.3 in~\cite{Wang2021some}.
\end{rem}

We arrange the rest of this paper as follows. In Section~\ref{Sec:top-pre}, we recall the notions of topological pressures and discuss their relations. Section~\ref{sec:mea-pre} presents various kinds of measure-theoretical pressures and proofs of Theorem~\ref{thm:vp-among-mp} and Theorem~\ref{thm:Billingsley-theorem}. We then prove Theorems~\ref{thm:variational-p-for-BP} and~\ref{thm:variational-p-for-P} in Section~\ref{sec:VP} and Theorem~\ref{thm:va-pp-gene} in Section~\ref{sec:va-of-gene}.

\section{Topological pressures}\label{Sec:top-pre}
In this section we recall four types of topological pressures: Pesin-Pitskel topological pressure, lower capacity topological pressure, upper capacity topological pressure, and packing topological pressure, and present some basic properties for these pressures.

For any $n\in\mathbb{N}$ and $x,y\in X$, let
\begin{equation*}\label{eq2.1}
d_n(x,y)=\max\{d(T^i(x),T^i(y)):0\leq i<n\}.
\end{equation*}
For $x\in X$ and $f\in C(X,\R)$, we write
\[f_n(x)=\sum_{i=0}^{n-1}f(T^i(x)).\]

For any $n\in\N$, $\varepsilon>0$, and $x\in X$, let $B_n(x,\varepsilon)=\{y\in X: d_n(x,y)<\varepsilon\}$ and $\overline{B}_{n}(x,\varepsilon)=\{y\in X:d_{n}(x,y)\leq\vp\}$. Given $f\in C(X,\R)$, let
\[f_n(x,\varepsilon)=\sup_{y\in B_n(x,\vp)}f_n(y),~~\text{and}~~\overline{f}_n(x,\varepsilon)=\sup_{y\in \overline{B}_n(x,\vp)}f_n(y).\]

Let $Z\subset X$ be a nonempty set. Given $n\in\N$, $\alpha\in\R$, $\varepsilon>0$, and $f\in C(X,\R)$, define
\begin{equation}\label{eq:definiton-of-sum-for-BP}
M(n,\alpha,\varepsilon,Z,T,f)=\inf\{\sum_{i}e^{-\alpha n_i+f_{n_i}(x_i)}: Z\subset\cup_{i}B_{n_i}(x_i,\varepsilon)\},
\end{equation}
where the infimum is taken over all finite or countable collections of $\{B_{n_i}(x_i,\varepsilon)\}_i$ such that $x_i\in X$, $n_i\geq n$ and $\bigcup_iB_{n_i}(x_i,\varepsilon)\supset Z$.
Likewise, we define
\begin{equation}\label{eq:definiton-of-sum-for-LUP}
R(n,\alpha,\varepsilon,Z,T,f)=\inf\{\sum_{i}e^{-\alpha n+f_{n}(x_i)}: Z\subset\cup_{i}B_{n}(x_i,\varepsilon)\},
\end{equation}
where the infimum is taken over all finite or countable collections of $\{B_{n}(x_i,\varepsilon)\}_i$ such that $x_i\in X$ and $\bigcup_iB_{n}(x_i,\varepsilon)\supset Z$.

Define
\begin{equation}\label{eq:definiton-of-sum-for-pp}
M^P(n,\alpha,\varepsilon,Z,T,f)=\sup\{\sum_{i}e^{-s n_i+f_{n_i}(x_i)}\},
\end{equation}
where the supremum is taken over all finite or countable pairwise disjoint families $\{\overline{B}_{n_i}(x_i,\varepsilon)\}$ such that $x_i\in Z$, $n_i\geq n$  for all $i$, where $\overline{B}_{n_i}(x_i,\varepsilon)=\{y\in X:d_{n_i}(x,y)\leq\vp\}$.

Let
\begin{align*}
M(\alpha,\varepsilon,Z,T,f)&=\lim_{n\to\infty}M(n,\alpha,\varepsilon,Z,T,f),\\
\underline{R}(\alpha,\varepsilon,Z,T,f)&=\liminf_{n\to\infty}R(n,\alpha,\varepsilon,Z,T,f),\\
\overline{R}(\alpha,\varepsilon,Z,T,f)&=\limsup_{n\to\infty}R(n,\alpha,\varepsilon,Z,T,f),\\
M^P(\alpha,\varepsilon,Z,T,f)&=\lim_{n\to\infty}M^P(n,\alpha,\varepsilon,Z,T,f).
\end{align*}
Define
\[M^{\calP}(\alpha,\varepsilon,Z,T,f)=\inf\{\sum_{i=1}^\infty M^P(\alpha,\varepsilon,Z_i,T,f):Z\subset\cup_{i=1}^\infty Z_i\}.\]

It is routine to check that when $\alpha$ goes from $-\infty$ to $+\infty$, the quantities
\[M(\alpha,\varepsilon,Z,T,f),~\underline{M}(\alpha,\varepsilon,Z,T,f),
\overline{M}(\alpha,\varepsilon,Z,T,f),~M^{\calP}(\alpha,\varepsilon,Z,T,f)\]
jump from $+\infty$ to $0$ at unique critical values respectively. Hence we can define the numbers
\begin{align*}
P^B(\varepsilon,Z,T,f)&=\sup\{\alpha:M(\alpha,\varepsilon,Z,T,f)=+\infty\}\\
&=\inf\{\alpha:M(\alpha,\varepsilon,Z,T,f)=0\},\\
\underline{CP}(\varepsilon,Z,T,f)&=\sup\{\alpha:\underline{R}(\alpha,\varepsilon,Z,T,f)=+\infty\}\\
&=\inf\{\alpha:\underline{R}(\alpha,\varepsilon,Z,T,f)=0\},\\
\overline{CP}(\varepsilon,Z,T,f)&=\sup\{\alpha:\overline{R}(\alpha,\varepsilon,Z,T,f)=+\infty\}\\
&=\inf\{\alpha:\overline{R}(\alpha,\varepsilon,Z,T,f)=0\},\\
P^P(\varepsilon,Z,T,f)&=\sup\{\alpha:M^{\calP}(\alpha,\varepsilon,Z,T,f)=+\infty\}\\
&=\inf\{\alpha:M^{\calP}(\alpha,\varepsilon,Z,T,f)=0\}.
\end{align*}

\begin{defn}\label{defn:pressures}
We call the following quantities
\begin{align*}
P^B(Z,T,f)&=\lim_{\vp\to0}P^B(\varepsilon,Z,T,f),\\
\underline{CP}(Z,T,f)&=\lim_{\vp\to0}\underline{CP}(\varepsilon,Z,T,f),\\
\overline{CP}(Z,T,f)&=\lim_{\vp\to0}\overline{CP}(\varepsilon,Z,T,f),\\
P^P(Z,T,f)&=\lim_{\vp\to0}P^P(\varepsilon,Z,T,f)
\end{align*}
\emph{Pesin-Pitskel}, \emph{lower capacity}, \emph{upper capacity}, and \emph{packing topological pressures of $T$ on the set $Z$ with respect to $f$}.
\end{defn}

\begin{rem}
The definitions of Pesin-Pitskel, lower capacity, and upper capacity topological pressures follow the generalized Carath\'{e}odory construction described in~\cite{Pesin1997Dimension}. For more details, see~\cite[p.74]{Pesin1997Dimension}. Wang and Chen in~\cite{Wang2012Variational} introduced packing topological pressure.
\end{rem}

Replacing $f_{n_i}(x_i)$ in Eqs.~\eqref{eq:definiton-of-sum-for-BP}, \eqref{eq:definiton-of-sum-for-pp} by $f_{n_i}(x_i,\vp)$ and $\overline{f}_{n_i}(x_i,\vp)$ respectively and $f_n(x_i)$ in Eq.\eqref{eq:definiton-of-sum-for-LUP} by $f_n(x_i,\vp)$, we can define new functions $\calM$, $\calR$, $\calM^P$. For any set $Z\subset X$ and $\vp>0$, we denote the respective critical values by
\[{P^B}'(\varepsilon,Z,T,f),~\underline{CP}'(\varepsilon,Z,T,f),
~\overline{CP}'(\varepsilon,Z,T,f),~{P^P}'(\varepsilon,Z,T,f).\]

\begin{prop}\label{prop:alternative-definition-of-BP}
Let $(X,T)$ be a TDS, $f\in C(X,\R)$, and $Z\subset X$. Then
\begin{align*}
P^B(Z,T,f)&
=\lim_{\vp\to0}{P^B}'(\varepsilon,Z,T,f),\\
\underline{CP}(Z,T,f)&
=\lim_{\vp\to0}\underline{CP}'(\varepsilon,Z,T,f),\\
\overline{CP}(Z,T,f)&
=\lim_{\vp\to0}\overline{CP}'(\varepsilon,Z,T,f),\\
P^P(Z,T,f)&=\lim_{\vp\to0}{P^P}'(\varepsilon,Z,T,f).
\end{align*}
\end{prop}

\begin{proof}
Fix $\vp>0$. It is clear that $P^B(\vp,Z,T,f)\leq{P^B}'(\varepsilon,Z,T,f)$. Let
\[\gamma(\vp)=\sup\{|f(x)-f(y)|:d(x,y)<2\vp\}.\]
A brief computation shows that
\[f_n(x,\vp)-n\gamma(\vp)\leq f_n(x),~\forall~n\in\N.\]
It then follows that
\begin{align*}
M(n,\alpha,\varepsilon,Z,T,f)&=\inf\{\sum_{i}e^{-\alpha n_i+f_{n_i}(x_i)}: Z\subset\cup_{i}B_{n_i}(x_i,\varepsilon)\}\\
&\geq\inf\{\sum_{i}e^{-\alpha n_i+f_{n_i}(x_i,\vp)-n_i\gamma(\vp)}: Z\subset\cup_{i}B_{n_i}(x_i,\varepsilon)\}\\
&=\inf\{\sum_{i}e^{-(\alpha+\gamma(\vp)) n_i+f_{n_i}(x_i,\vp)}: Z\subset\cup_{i}B_{n_i}(x_i,\varepsilon)\}\\
&=\calM(n,\alpha+\gamma(\vp),\varepsilon,Z,T,f).
\end{align*}
Letting $n\to\infty$ yields
\[M(\alpha,\varepsilon,Z,T,f)\geq \calM(\alpha+\gamma(\vp),\varepsilon,Z,T,f).\]
This implies that
\[P^B(\varepsilon,Z,T,f)\geq {P^B}'(\varepsilon,Z,T,f)-\gamma(\vp).\]
It then follows that
\[P^B(\vp,Z,T,f)\leq{P^B}'(\varepsilon,Z,T,f)\leq P^B(\varepsilon,Z,T,f)+\gamma(\vp). \]
Since $f$ is uniformly continuous on $X$, the desired equality follows by letting $\vp\to0$.

The other equalities can be proven similarly.
\end{proof}

The following are some basic properties of these pressures.

\begin{prop}\label{prop:properties-of-pressures}
Let $(X,T)$ be a TDS and $f\in C(X,\R)$. Then the following assertions hold:
\begin{itemize}
  \item[1.] For any $Z\subset X$, $P^B(Z,T,f)\leq \underline{CP}(Z,T,f)\leq \overline{CP}(Z,T,f)$.
  \item[2.] If $Z_1\subset Z_2$, then $\scrP(Z_1,f)\leq\scrP(Z_2,f)$, where $\scrP\in\{P^B,\underline{CP},\overline{CP},P^P\}$.
  \item[3.] If $Z=\cup_{i\in I}Z_i$ is a union of sets $Z_i\subset X$, with $I$ at most countable, then
      \begin{itemize}
        \item[3-a.]\label{thm:sub-add-for-P} $M(\alpha,\varepsilon,Z,T,f)\leq\sum_{i}M(\alpha,\varepsilon,Z_i,T,f)$;
        \item[3-b.] $M^{\calP}(\alpha,\varepsilon,Z,T,f)\leq\sum_{i}M^{\calP}
        (\alpha,\varepsilon,Z_i,T,f)$;
        \item[3-c.] $P^B(Z,T,f)=\sup_{i\in I}P^B(Z_i,T,f)$;
        \item[3-d.] $P^P(Z,T,f)=\sup_{i\in I}P^P(Z_i,T,f)$;
        \item[3-e.] $\underline{CP}(Z,T,f)\geq\sup_{i\in I}\underline{CP}(Z_i,T,f)$;
        \item[3-f.] $\overline{CP}(Z,T,f)\geq\sup_{i\in I}\overline{CP}(Z_i,T,f)$.
      \end{itemize}
  \item[4.]\label{thm:p-leq-pp-leq-up} For any $Z\subset X$, $P^B(Z,T,f)\leq P^P(Z,T,f)\leq \overline{CP}(Z,T,f)$.
  \item[5.] If $Z$ is $T$-invariant and compact, then
  \[P^B(Z,T,f)=P^P(Z,T,f)=\underline{CP}(Z,T,f)=\overline{CP}(Z,T,f).\]
\end{itemize}
\end{prop}

\begin{proof}
1 and 2 follow directly from the definitions of topological pressures. 3-a follows from (3) of Proposition 1 in~\cite{Pesin1997Dimension}. 3-c,3-e, and 3-f follows from Theorems 11.2 and 11.3 in~\cite{Pesin1997Dimension}.

We now show 3-b. Given $\gamma>0$ and $i\in I$, we can find $\{Z_{i,j}\}_{j\geq0}$ such that $Z_i\subset\cup_{j\geq0}Z_{i,j}$ and
\[\sum_{j\geq0} M^P(\alpha,\varepsilon,Z_{i,j},T,f)\leq M^{\calP}(\alpha,\varepsilon,Z_i,T,f)+\frac{\gamma}{2^i}.\]
Taking sums over all $i\in I$ yields
\[\sum_{i\in I}\sum_{j\geq0} M^P(\alpha,\varepsilon,Z_{i,j},T,f)\leq \sum_{i\in I}M^{\calP}(\alpha,\varepsilon,Z_i,T,f)+3\gamma.\]
The desired inequality then follows.

We now show that 3-d holds. If $\sup_{i\in I}P^P(Z_i,T,f)<s$ then for any $\vp>0$, $P^P(\vp,Z_i,T,f)<s$. Thus $M^{\calP}(s,\vp,Z_i,T,f)=0$ and we obtaind $M^{\calP}(s,\vp,Z,T,f)=0$ by utilizing 3-b. Hence $P^P(\vp,Z,T,f)\leq s$. It follows that $P^P(\vp,Z,T,f)\leq \sup_{i\in I}P^P(Z_i,T,f)$. The opposite inequality follows from 2.

We shall modify slightly the argument of \cite[Proposition 2.1]{feng2012variational} to prove 4.
We first show that $P^B(Z,T,f)\leq P^P(Z,T,f)$. Suppose that $P^B(Z,T,f)>s>-\infty$. For any $\vp>0$ and $n\in\N$, let
\[\F_{n,\vp}=\{\calF:\calF=\{\overline{B}_n(x_i,\vp)\}, x_i\in Z,~\text{and}~\calF~\text{is a disjoint family}\}.\]
Take $\calF(n,\vp,Z)\in\F_{n,\vp}$ such that
$|\calF(n,\vp,Z)|=\max_{\calF\in\F_{n,\vp}}\{|\calF|\}$, where $|\calF|$ denotes the cardinality of $\calF$.
For convenience, we denote $\calF(n,\vp,Z)=\{\overline{B}_n(x_i,\vp),i=1,\ldots,|\calF(n,\vp,Z)|\}.$
It is easy to check that
\[Z\subset\bigcup_{i=1}^{|\calF(n,\vp,Z)|}B_n(x_i,2\vp+\delta),~\forall \delta>0.\]
Then for any $s\in\R$,
\[M(n,s,2\vp+\delta,Z,T,f)\leq e^{-ns}\sum_{i=1}^{|\calF(n,\vp,Z)|}f_n(x_i)\leq M^P(n,s,\vp,Z,T,f).\]
It thus follows that $M(s,2\vp+\delta,Z,T,f)\leq M^P(s,\vp,Z,T,f)$. By 3-a, we have $M(s,2\vp+\delta,Z,T,f)\leq M^{\calP}(s,\vp,Z,T,f)$.
Since $P^B(Z,T,f)>s>-\infty$, $M(s,2\vp+\delta,Z,T,f)\geq1$ when $\vp$ and $\delta$ are small enough. Thus $M^{\calP}(s,\vp,Z,T,f)\geq1$.
This implies that $P^{P}(\vp,Z,T,f)\geq s$ for $\vp$ small enough. Hence $P^{P}(Z,T,f)\geq s$ and $P^B(Z,T,f)\leq P^P(Z,T,f)$.

We shall show $P^P(Z,T,f)\leq \overline{CP}(Z,T,f)$.

Case 1. $P^P(Z,T,f)=-\infty$. It is obvious that the inequality holds.

Case 2.  $P^P(Z,T,f)>-\infty$. Choose $-\infty<t<s<P^P(Z,T,f)$. Then there exists $\delta>0$, such that for any $\vp\in(0,\delta)$, $P^P(\vp,Z,T,f)>s$ and $M^P(s,\vp,Z,T,f)\geq M^{\calP}(s,\vp,Z,T,f)=\infty$. Hence, for any $N\in\N$, there exists a countable pairwise disjoint family $\{\overline{B}_{n_i}(x_i,\vp)\}$ such that $x_i\in Z$, $n_i\geq N$ for all $i$, and $\sum_{i}e^{-n_is+f_{n_i}(x)}>1$. For each $k$, let
\[m_k=\{x_i:n_i=k\}.\]
Then
\[\sum_{k=N}^\infty \sum_{x\in m_k}e^{f_k(x)} e^{-k s}>1.\]
It is easy to check that there exists $k\geq N$ such that $\sum_{x\in m_k}e^{f_k(x)} e^{-k t}\geq 1-e^{t-s}$ (otherwise, $\sum_{k=N}^\infty \sum_{x\in m_k}e^{f_k(x)} e^{-k s}\leq 1$). Fixing a collection $\{B_k(y_i,\frac{\vp}{2})\}_{i\in I}$ such that $Z\subset \cup_{i\in I}B(y_i,\frac{\vp}{2})$, where $I$ is at most countable, it is not difficult to check that for any $x_1,x_2\in m_k$ there exit distinguish $y_1,y_2$ such that $x_i\in B(y_i,\frac{\vp}{2})$, $i=1,2$.
Then
\[\calR(k,t,\frac{\vp}{2},Z,T,f)\geq\sum_{x\in m_k}e^{f_k(x)} e^{-k t}\geq 1-e^{t-s}.\]
Hence
\[\overline{\calR}(t,\frac{\vp}{2},Z,T,f)=\limsup_{n\to\infty}\calR(k,t,\frac{\vp}{2},Z,T,f)\geq 1-e^{t-s}>0.\]
Thus $\overline{CP}'(\frac{\vp}{2},Z,T,f)\geq t$.
Letting  $\vp\to0$ yields $\overline{CP}(Z,T,f)\geq t$. Since $t\in(-\infty,P^P(Z,T,f))$, it follows that $P^P(Z,T,f)\leq \overline{CP}(Z,T,f)$.

(5) follows from (2) of Theorem 11.5 in~\cite{Pesin1997Dimension}, where $P^B(Z,T,f)=\underline{CP}(Z,T,f)=\overline{CP}(Z,T,f)$.
\end{proof}


\section{Measure-theoretic pressures}\label{sec:mea-pre}
In this section, we discuss the relations among various types of measure-theoretic topological pressure.

Let $(X,T)$ be a TDS, $f\in C(X,\R)$ and $\mu\in\mathcal{M}(X)$. The \emph{measure-theoretic lower and upper local pressures of $x\in X$ with respect to $\mu$ and $f$} are defined by
\begin{align*}
  \underline{P}_{\mu}(x,T,f)&:=\lim_{\vp\to0}\liminf_{n\to\infty}\frac{-\log\mu(B_n(x,\vp))+f_n(x)}{n},\\
  \overline{P}_{\mu}(x,T,f)&:=\lim_{\vp\to0}\limsup_{n\to\infty}\frac{-\log\mu(B_n(x,\vp))+f_n(x)}{n}.
\end{align*}

\begin{defn}\label{defn:M-L-U-pressures}
The \emph{measure-theoretic lower and upper local pressures of $\mu$ with respect to $f$} are defined as
\begin{align*}
\underline{P}_{\mu}(T,f):=\int \underline{P}_{\mu}(x,T,f)\md\mu(x),\\
\overline{P}_{\mu}(T,f):=\int \overline{P}_{\mu}(x,T,f)\md\mu(x).
\end{align*}
\end{defn}

\begin{defn}\label{defn:M-pressures}
We call the following quantities
\begin{align*}
P_\mu^{B}(T,f)&:=\lim_{\vp\to0}\lim_{\delta\to0}\inf\{P^B(\vp,Z,T,f):\mu(Z)\geq1-\delta\}\\
&~=\lim_{\vp\to0}\lim_{\delta\to0}\inf\{{P^B}'(\vp,Z,T,f):\mu(Z)\geq1-\delta\}\\
\underline{CP}_{\mu}(T,f)&:=\lim_{\vp\to0}\lim_{\delta\to0}\inf\{\underline{CP}(\vp,Z,T,f):\mu(Z)\geq1-\delta\}\\
&~=\lim_{\vp\to0}\lim_{\delta\to0}\inf\{\underline{CP}'(\vp,Z,T,f):\mu(Z)\geq1-\delta\}\\
\overline{CP}_{\mu}(T,f)&:=\lim_{\vp\to0}\lim_{\delta\to0}\inf\{\overline{CP}(\vp,Z,T,f):\mu(Z)\geq1-\delta\}\\
&~=\lim_{\vp\to0}\lim_{\delta\to0}\inf\{\overline{CP}'(\vp,Z,T,f):\mu(Z)\geq1-\delta\}\\
P_{\mu}^{P}(T,f)&:=\lim_{\vp\to0}\lim_{\delta\to0}\inf\{P^P(\vp,Z,T,f):\mu(Z)\geq1-\delta\}\\
&~=\lim_{\vp\to0}\lim_{\delta\to0}\inf\{{P^P}'(\vp,Z,T,f):\mu(Z)\geq1-\delta\}
\end{align*}
\emph{Pesin-Pitskel, lower capacity, upper capacity, and packing pressures of $\mu$ with respect to $f$}.
\end{defn}

Katok in~\cite{Katok1980Lyapunov} introduced a type of measure-theoretic entropy. Recently, Wang in~\cite{Wang2021some} studied the dimension types of this entropy. Following along the line of topological pressures in Section~\ref{Sec:top-pre}, we shall introduce four dimension types of measure-theoretic pressure in the sense of Katok.

Let $Z\subset X$ be a nonempty set. Given $\mu\in\calM(X)$, $n\in\N$, $\alpha\in\R$, $\varepsilon>0$, $0<\delta<1$, and $f\in C(X,\R)$, define
\begin{equation}\label{eq:definiton-of-sum-for-KBP}
M_\mu(n,\alpha,\varepsilon,\delta,T,f)=\inf\{\sum_{i}e^{-\alpha n_i+f_{n_i}(x_i)}: \mu\big(\cup_{i}B_{n_i}(x_i,\varepsilon)\big)\geq 1-\delta\},
\end{equation}
where the infimum is taken over all finite or countable collections of $\{B_{n_i}(x_i,\varepsilon)\}_i$ such that $x_i\in X$, $n_i\geq n$ and $\mu\big(\cup_{i}B_{n_i}(x_i,\varepsilon)\big)\geq 1-\delta$.
Likewise, we define
\begin{equation}\label{eq:definiton-of-sum-for-LUKP}
R_\mu(n,\alpha,\varepsilon,\delta,T,f)=\inf\{\sum_{i}e^{-\alpha n+f_{n}(x_i)}: \mu\big(\cup_{i}B_{n}(x_i,\varepsilon)\big)\geq1-\delta\},
\end{equation}
where the infimum is taken over all finite or countable collections of $\{B_{n}(x_i,\varepsilon)\}_i$ such that $x_i\in X$ and $\mu\big(\cup_{i}B_{n}(x_i,\varepsilon)\big)\geq1-\delta$.

Let
\begin{align*}
M_\mu(\alpha,\varepsilon,\delta,T,f)&=\lim_{n\to\infty}M_\mu(n,\alpha,\varepsilon,\delta,T,f),\\
\underline{M}_\mu(\alpha,\varepsilon,\delta,T,f)&=\liminf_{n\to\infty}R_\mu(n,\alpha,\varepsilon,\delta,T,f),\\
\overline{M}_\mu(\alpha,\varepsilon,\delta,T,f)&=\limsup_{n\to\infty}R_\mu(n,\alpha,\varepsilon,\delta,T,f).
\end{align*}
Define
\[M^{\calP}_\mu(\alpha,\varepsilon,\delta,T,f)=\inf\{\sum_{i=1}^\infty M^P(\alpha,\varepsilon,Z_i,T,f):\mu\big(\cup_{i=1}^\infty Z_i\big)\geq1-\delta\}.\]

Thus, when $\alpha$ goes from $-\infty$ to $+\infty$, the quantities
\[M_\mu(\alpha,\varepsilon,\delta,T,f),~\underline{M}_\mu(\alpha,\varepsilon,\delta,T,f),
\overline{M}_\mu(\alpha,\varepsilon,\delta,T,f),\text{and}~M^{\calP}_\mu(\alpha,\varepsilon,\delta,T,f)\]
jump from $+\infty$ to $0$ at unique critical values respectively. Hence we can define the numbers
\begin{align*}
P^{KB}_\mu(\varepsilon,\delta,T,f)&=\sup\{\alpha:M_\mu(\alpha,\varepsilon,\delta,T,f)=+\infty\}\\
&=\inf\{\alpha:M_\mu(\alpha,\varepsilon,\delta,T,f)=0\},\\
\underline{CP}^K_\mu(\varepsilon,\delta,T,f)&=\sup\{\alpha:\underline{R}_\mu(\alpha,\varepsilon,\delta,T,f)=+\infty\}\\
&=\inf\{\alpha:\underline{R}_\mu(\alpha,\varepsilon,\delta,T,f)=0\},\\
\overline{CP}^K_\mu(\varepsilon,\delta,T,f)&=\sup\{\alpha:\overline{R}_\mu(\alpha,\varepsilon,\delta,T,f)=+\infty\}\\
&=\inf\{\alpha:\overline{R}_\mu(\alpha,\varepsilon,\delta,T,f)=0\},\\
P^{KP}_\mu(\varepsilon,\delta,T,f)&=\sup\{\alpha:M^{\calP}_\mu(\alpha,\varepsilon,\delta,T,f)=+\infty\}\\
&=\inf\{\alpha:M^{\calP}_\mu(\alpha,\varepsilon,\delta,T,f)=0\}.
\end{align*}

\begin{defn}\label{defn:M-K-pressures}
We call the following quantities
\begin{align*}
P^{KB}_\mu(T,f)&=\lim_{\vp\to0}\lim_{\delta\to0}P^{KB}_\mu(\varepsilon,\delta,T,f),\\
\underline{CP}^K_\mu(T,f)&=\lim_{\vp\to0}\lim_{\delta\to0}\underline{CP}^K_\mu(\varepsilon,\delta,T,f),\\
\overline{CP}^K_\mu(T,f)&=\lim_{\vp\to0}\lim_{\delta\to0}\overline{CP}^K_\mu(\varepsilon,\delta,T,f),\\
P^{KP}_\mu(T,f)&=\lim_{\vp\to0}\lim_{\delta\to0}P^{KP}_\mu(\varepsilon,\delta,T,f)
\end{align*}
\emph{Pesin-Pitskel}, \emph{lower capacity}, \emph{upper capacity}, and \emph{packing  pressures of $\mu$ in the sense of Katok with respect to $f$}, respectively.
\end{defn}

\begin{rem}
If we replace $f_{n_i}(x_i)$ in Eqs.~\eqref{eq:definiton-of-sum-for-KBP}, \eqref{eq:definiton-of-sum-for-pp} by $f_{n_i}(x_i,\vp)$ and $\overline{f}_{n_i}(x_i,\vp)$ respectively and $f_n(x_i,\vp)$ by $f_n(x_i)$ in Eq.\eqref{eq:definiton-of-sum-for-LUKP}, we can define new functions $\calM_\mu$, $\calR_\mu$, $\calM^{\calP}_\mu$. For any $\vp>0$ and $0<\delta<1$, we denote the respective critical values by
\[{P^{KB}_\mu}'(\varepsilon,\delta,T,f),~{\underline{CP}^K_\mu}'(\varepsilon,\delta,T,f),
~{\overline{CP}^K_\mu}'(\varepsilon,\delta,T,f),~{P^{KP}_\mu}'(\varepsilon,\delta,T,f).\]
\end{rem}

\begin{prop}\label{prop:alternative-definition-of-KBP}
Let $(X,T)$ be a TDS, $\mu\in \calM(X)$, and $f\in C(X,\R)$. Then
\begin{align*}
P^{KB}_\mu(T,f)&
=\lim_{\vp\to0}\lim_{\delta\to0}{P^{KB}_\mu}'(\varepsilon,\delta,T,f),\\
\underline{CP}^K_\mu(T,f)&=\lim_{\vp\to0}\lim_{\delta\to0}{\underline{CP}^K_\mu}'(\varepsilon,\delta,T,f),\\
\overline{CP}^K_\mu(T,f)&=\lim_{\vp\to0}\lim_{\delta\to0}{\overline{CP}^K_\mu}'(\varepsilon,\delta,T,f),\\
P^{KP}_\mu(T,f)&=\lim_{\vp\to0}\lim_{\delta\to0}{P^{KP}_\mu}'(\varepsilon,\delta,T,f)
\end{align*}
\end{prop}

\begin{proof}
The proof is analogous to that of Proposition~\ref{prop:alternative-definition-of-BP}, so we omit it.
\end{proof}

\begin{prop}\label{prop:MKPs-equal-MPs}
Let $(X,T)$ be a TDS, $\mu\in\calM(X)$, and $f\in C(X,\R)$. Then
\begin{align*}
P^{KB}_\mu(T,f)&= P^{B}_\mu(T,f),~\underline{CP}^K_\mu(T,f)= \underline{CP}_\mu(T,f),\\
\overline{CP}^K_\mu(T,f)&\leq\overline{CP}_\mu(T,f),~P^{KP}_\mu(T,f)=P^{P}_\mu(T,f).
\end{align*}
\end{prop}

\begin{proof}
We shall show that $P^{KB}_\mu(T,f)\leq P^{B}_\mu(T,f)$. For any $n\in\N$, $\alpha\in\R$, $\vp>0$, $0<\delta<1$, and $Z$ with $\mu(Z)\geq 1-\delta$,
\[M_\mu(n,\alpha,\varepsilon,\delta,T,f)\leq M(n,\alpha,\varepsilon,Z,T,f).\]
Letting $n\to\infty$ yields
\[M_\mu(\alpha,\varepsilon,\delta,T,f)\leq M(\alpha,\varepsilon,Z,T,f).\]
This shows that
\[P^{KB}_\mu(\varepsilon,\delta,T,f)\leq P^B(\vp,Z,T,f),\]
and consequently
\[P^{KB}_\mu(\varepsilon,\delta,T,f)\leq\inf\{P^B(\vp,Z,T,f):\mu(Z)\geq1-\delta\}.\]
Letting $\delta\to0$ and $\vp\to0$, the desired inequality follows.

We can prove similarly $\underline{CP}^K_\mu(T,f)\leq\underline{CP}_\mu(T,f)$ and $
\overline{CP}^K_\mu(T,f)\leq\overline{CP}_\mu(T,f)$.

To prove $P^{KB}_\mu(T,f)\geq P^{B}_\mu(T,f)$, let $a=P^{KB}_\mu(T,f)$. For any $s>0$, there exists $\vp'>0$ such that
\[\lim_{\delta\to0}P^{KB}_\mu(\varepsilon,\delta,T,f)<a+s,~\forall~\vp<\vp'.\]
It follows that for any $\vp\in(0,\vp')$, there exists $\delta_\vp$ so that
\[P^{KB}_\mu(\varepsilon,\delta,T,f)<a+s,~\forall~\delta<\delta_\vp.\]
This implies that $\lim_{n\to\infty}M_\mu(n,a+s,\varepsilon,\delta,T,f)=0$. For any $N\in\N$, we can find a sequence of $\delta_{N,m}$ with $\lim_{m\to0}\delta_{N,m}=0$ and a collection of $\{B_{n_i}(x_i,\varepsilon)\}_{i\in I_{N,m}}$ such that $x_i\in X$, $n_i\geq N$, $\mu\big(\cup_{i\in I_{N,m}}B_{n_i}(x_i,\varepsilon)\big)\geq 1-\delta_{N,m}$, and
\[\sum_{i\in I_{N,m}}e^{-(a+s) n_i+f_{n_i}(x_i)}\leq \frac{1}{2^m}.\]
Let
\[Z_N=\bigcup_{m\in\N}\bigcup_{i\in I_{N,m}}B_{n_i}(x_i,\varepsilon).\]
Then $\mu(Z_N)=1$ and
\[M(N,a+s,\varepsilon,Z_N,T,f)\leq 1.\]
Let $Z_\vp=\cap_{N\in\N}Z_N$. Thus $\mu(Z_\vp)=1$ and
\[M(N,a+s,\varepsilon,Z_\vp,T,f)\leq M(N,a+s,\varepsilon,Z_N,T,f)\leq 1,~\forall N\in\N.\]
It follows that
\[P^B(\vp,Z_\vp,T,f)\leq a+s.\]
Therefore,
\[P_\mu^{B}(T,f)=\lim_{\vp\to0}\lim_{\delta\to0}\inf\{P^B(\vp,Z,T,f):\mu(Z)\geq1-\delta\}\leq a+s.\]
The arbitrariness of $s$ then implies the desired inequality.

To prove $\underline{CP}^K_\mu(T,f)\geq\underline{CP}_\mu(T,f)$, let $a=\underline{CP}^K_\mu(T,f)$. For any $s>0$, there exists $\vp'>0$ such that for any $\vp\in(0,\vp')$, there exists $\delta_\vp$ so that
\[\liminf_{n\to\infty}R_\mu(n,a+s,\varepsilon,\delta,T,f)=0,~\forall~\delta<\delta_\vp.\]
Fix $\delta\in(0,\delta_\vp)$. For any $m\in\N$, we have
\[\liminf_{n\to\infty}R_\mu(n,a+s,\varepsilon,\frac{\delta}{2^m},T,f)=0.\]
Then for every $m\in\N$ there exits a family $\{B_{k_m}(x_i,\vp)\}_{i\in I_m}$ with $\mu(\cup_{i\in I_m}B_{k_m}(x_i,\vp))\geq1-\frac{\delta}{2^m}$ such that
\[\sum_{i\in I_m}e^{-(a+s) k_m+f_{k_m}(x_i)}\leq 1.\]
Let $Z_\delta=\cap_{m\in\N}\cup_{i\in I_m}B_{k_m}(x_i,\vp)$. Then $\mu(Z_\delta)\geq 1-\delta$. It is easy to check that
\[\liminf_{n\to\infty} R(n,a+s,\vp,Z_\delta,T,f)\leq 1.\]
Thus
\[\underline{CP}(\vp,Z_\delta,T,f)\leq a+s.\]
This implies that $\underline{CP}_\mu(T,f)\leq a+s$, and the desired inequality follows from the arbitrariness of $s$.

We now show the fourth equality. We first prove that $P^{P}_\mu(T,f)\geq P^{KP}_\mu(T,f)$. For any $s<P^{KP}_\mu(T,f)$, there exit $\vp',\delta'>0$ such that
\[P^{KP}_\mu(\varepsilon,\delta,T,f)>s,~\forall~\vp\in(0,\vp'),~\delta\in(0,\delta').\]
Thus
\[M^{\calP}_\mu(s,\varepsilon,\delta,T,f)=\infty.\]
For any $Z$ with $\mu(Z)\geq1-\delta$. If $Z\subset\cup_iZ_i$ then $\mu(\cup_iZ_i)\geq1-\delta$. It follows that
\[\sum_{i=1}^\infty M^P(s,\varepsilon,Z_i,T,f)=\infty,\]
which implies that $M^\calP(s,\varepsilon,Z,T,f)=\infty$. Hence $P^P(\vp,Z,T,f)\geq s$ and $P^{P}_\mu(T,f)\geq s$. This shows that $P^{P}_\mu(T,f)\geq P^{KP}_\mu(T,f)$.

We shall show the inverse inequality. If $s<P^{P}_\mu(T,f)$ then there exit $\vp',\delta'>0$ such that
\[\inf\{P^P(\varepsilon,Z,T,f):\mu(Z)\geq 1-\delta\}>s,~\forall~\vp\in(0,\vp'),~\delta\in(0,\delta').\]
For any family $\{Z_i\}_{i\geq1}$ with $\mu\big(\cup_iZ_i\big)\geq 1-\delta$, we have
\[P^P(\varepsilon,\cup_iZ_i,T,f)>s.\]
This implies that
\[M^\calP(s,\varepsilon,\cup_iZ_i,T,f)=\infty.\]
Thus
\[\sum_iM^P(s,\varepsilon,Z_i,T,f)=\infty.\]
Then
\[M^{\calP}_\mu(s,\varepsilon,\delta,T,f)=\infty.\]
Hence
\[P^{KP}_\mu(\varepsilon,\delta,T,f)>s,\]
which yields the desired inequality.
\end{proof}

\begin{rem}
Theorem~\ref{thm:vp-among-mp} states that $\overline{CP}^K_\mu(T,f)=\overline{CP}_\mu(T,f)$ for $\mu\in\calE(X,T)$. But it is not clear if $\overline{CP}^K_\mu(T,f)=\overline{CP}_\mu(T,f)$ for any $\mu\in\calM(X)$?
\end{rem}

\textbf{Proof of Theorem~\ref{thm:vp-among-mp}.}
Employing Theorem 11.6 in~\cite{Pesin1997Dimension}, the following equalities follow
\begin{equation}\label{eq:star}
P_\mu^{B}(T,f)=\underline{CP}_{\mu}(T,f)=\overline{CP}_{\mu}(T,f)=h_{\mu}(T)+\int f\md\mu.\tag{$\ast$}
\end{equation}
From the proof of 4 of Proposition~\ref{prop:properties-of-pressures}, it is easy to check that for any $Z\subset X$ and $\vp>0$,
\[P^B(3\vp,Z,T,f)\leq P^P(\vp,Z,T,f)\leq \overline{CP}(\vp,Z,T,f).\]
Thus
\[P_\mu^{B}(T,f)\leq P_{\mu}^{P}(T,f)\leq\overline{CP}_{\mu}(T,f),\]
which together with \eqref{eq:star} yields
\[P_{\mu}^{P}(T,f)=h_{\mu}(T)+\int f\md\mu.\]
Using Proposition~\ref{prop:MKPs-equal-MPs}, we obtain
\[P_{\mu}^{KP}(T,f)=P_{\mu}^{P}(T,f)=h_{\mu}(T)+\int f\md\mu\]
and
\[P_\mu^{B}(T,f)= P_\mu^{KB}(T,f)\leq\underline{CP}^K_{\mu}(T,f)\leq\overline{CP}^K_{\mu}(T,f)\leq \overline{CP}_{\mu}(T,f).\]
The equalities in Theorem~\ref{thm:vp-among-mp} then follows.
\hfill$\Box$

Before the proof of Theorem~\ref{thm:Billingsley-theorem}, we need the following lemma.
\begin{lem}(\cite[Theorem 2.1]{Mattila1995Geometry})\label{lem:cover-5r-lem}
Let $(X,d)$ be a compact metric space and $\calB=\{B(x_i, r_i)\}_{i\in \calI}$ be a family of closed (or open) balls in $X$. Then there exists a finite or countable subfamily $\calB'=\{B(x_i, r_i)\}_{i\in\calI'}$ of pairwise disjoint balls in $\calB$ such that
\[\bigcup_{B\in\calB}B\subset\bigcup_{i\in\calI'}B(x_i,5r_i).\]
\end{lem}

\textbf{Proof of Theorem~\ref{thm:Billingsley-theorem}.}
We now prove the first assertion. For a fixed $\beta>s$, let
\[Z_m=\{x\in Z:\limsup_{n\to\infty}\frac{-\log\mu(B_n(x,\vp))+f_n(x)}{n}\leq \frac{\beta+s}{2},~\forall~\vp\in(0,\frac1{m})\}.\]
It is easy to check that $Z=\cup_{m=1}^\infty Z_m$ using the fact that $\overline{P}_{\mu}(x,T,f)\leq s$ for all $x\in Z$. Fix $m\geq1$ and $\vp\in(0,\frac1{m})$. For each $x\in Z_m$, there exists $N\in\N$ such that
\[\mu\big(B_n(x,\vp)\big)\geq e^{-n(\frac{\beta+s}{2})+f_{n}(x)},~\forall~n\geq N.\]
Let
\[Z_{m,N}=\{x\in Z_m:\mu\big(B_n(x,\vp)\big)\geq e^{[-n(\beta+s)]/2+f_{n}(x)},~\forall~n\geq N\}.\]
It is clear that $Z_m=\cup_{N=1}^\infty Z_{m,N}$. Fix $N$ and $L\geq N$. For a finite or countable disjoint family $\calF=\{\overline{B}_{n_i}(x_i,\vp)\}_{i}$, where $x_i\in Z_{m,N}$, $n_i\geq L$,
\begin{align*}
\sum_{i}e^{-n_i\beta+f_{n_i}(x_i)}
&=\sum_{i}e^{-n_i[(\beta+s)/2+(\beta-s)/2]+f_{n_i}(x_i)}\\
&\leq e^{-L[(\beta-s)/2]}\sum_{i}e^{-n_i[(\beta+s)/2]+f_{n_i}(x_i)}\\
&\leq e^{-L[(\beta-s)/2]}\sum_{i}\mu(\overline{B}_{n_i}(x_i,\vp))\\
&\leq e^{-L[(\beta-s)/2]}.
\end{align*}
Since $\calF$ is arbitrary, it follows that
\[M^P(L,\beta,\vp,Z_{m,N},T,f)\leq e^{-L[(\beta-s)/2]}.\]
Hence $M^P(\beta,\vp,Z_{m,N},f)=0$. This implies that $M^{\calP}(\beta,\vp,Z_{m},T,f)=0$. It follows that $P^P(Z_m,f)\leq\beta$. Thus
\[P^P(Z,T,f)=\sup_mP^P(Z_m,T,f)\leq\beta.\] By the arbitrariness of $\beta$, we have $P^P(Z,T,f)\leq s$.

We shall show the second assertion. Fix $\beta<s$. Let $\delta=\frac{s-\beta}{2}$. Since ~$\overline{P}_\mu(x,T,f)\geq s$, there exists $\vp>0$ such that
\[\limsup_{n\to\infty}\frac{-\log\mu(B_n(x,\vp))+f_n(x)}{n}>\beta+\delta.\]
Next we show that $M^{\calP}(s,\frac{\vp}{5},Z,T,f)=\infty$ which implies that $P^P(Z,T,f)\geq P^P(\frac{\vp}{5},Z,T,f)\geq s$. To this end, it suffices to show that $M^P(s,\vp/5,E,f)=\infty$ for any Borel $E\subset Z$ with $\mu(E)>0$. Fix such a set $E$. Let
\[E_n=\{x\in E:\mu(B_n^\alpha(x,\vp))<e^{-n(\beta+\delta)+f_n(x)}\},~~n\in\N.\]
It is clear that $E=\cup_{n=N}^\infty E_n$ for each $N\in\N$. Fix $N\in\N$. Then $\mu(\cup_{n=N}^\infty E_n)=\mu(E)$. Hence there exists $n\geq N$ such tat
\[\mu(E_n)\geq\frac{1}{n(n+1)}\mu(E).\]
Fix such $n$ and let $\calB=\{B_n(x,\frac{\vp}{5}):x\in E_n\}$. By Lemma~\ref{lem:cover-5r-lem} (in which we use $d_n$ instead of $d$), there exists a finite or countable pairwise disjoint family $\{B_n(x_i,\frac{\vp}{5})\}$ with $x_i\in E_n$ such that
\[E_n\subset\bigcup_{x\in E_n}B_n(x,\frac{\vp}{5})\subset\bigcup_iB_n(x_i,\vp).\]
Hence
\begin{align*}
M^P(N,\beta,\frac{\vp}{5},E,T,f)&\geq M^P(N,\beta,\frac{\vp}{5},E_n,T,f)\geq\sum_{i}e^{-n\beta+f_n(x_i)}\geq e^{n\delta}\sum_{i}e^{-n(\beta+\delta)+f_n(x_i)}\\
&\geq e^{n\delta}\sum_{i}\mu(B_n(x_i,\vp))\geq e^{n\delta}\mu(E_n)\geq\frac{e^{n\delta}}{n(n+1)}\mu(E_n).
\end{align*}
Since $\frac{e^{n\delta}}{n(n+1)}\to\infty$ as $n\to\infty$, it follows by letting $N\to\infty$ that $M^P(\beta,\frac{\vp}{5},E,T,f)=\infty$.
\hfill$\Box$

\begin{prop}\label{prop:m-p-leq-kp}
Let $(X,T)$ be a TDS and $f\in C(X,\R)$. For any $\mu\in\calM(X)$, we have (1) $\underline{P}_u(T,f)\leq P_\mu^{KB}(T,f)$;
(2) $\overline{P}_u(T,f)\leq P_{\mu}^{KP}(T,f)$.
\end{prop}

\begin{proof}
We now show the first inequality. Fix $s<\underline{P}_u(T,f)$. Then there exists a Borel set $Z\subset X$ with $\mu(Z)>0$ such that
\[\lim_{\vp\to0}\liminf_{n\to\infty}\frac{-\log\mu(B_n(x,\vp))+f_n(x)}{n}>s.\]
For each $m\geq1$, put
\[Z_m=\{x\in Z:\liminf_{n\to\infty}\frac{-\log\mu(B_n(x,\vp))+f_n(x)}{n}>s,~\forall~\vp\in(0,\frac1{m}]\}.\]
Since $\frac{-\log\mu(B_n(x,\vp))+f_n(x)}{n}$ increases when
$\varepsilon$ decreases, it follows that
\[Z_m=\{x\in Z:\liminf_{n\to\infty}\frac{-\log\mu(B_n(x,\vp))+f_n(x)}{n}>s,~\vp=\frac1{m}\}.\]
Then $Z_m\subset Z_{m+1}$ and $\cup_{m=1}^\infty Z_m=Z$. Using the continuity of measures yields
\[\lim_{m\to\infty} \mu(Z_m)=\mu(Z).\]
Take $M\geq1$ with $\mu(Z_M)>\frac1{2}\mu(Z)$. For every $N\geq1$, put
\begin{align*}
Z_{M,N}=&\{x\in Z_M:\frac{-\log\mu(B_n(x,\vp))+f_n(x)}{n}>s,~\forall~n\geq N,\vp\in(0,\frac1{M}]\}\\
=&\{x\in Z_M:\frac{-\log\mu(B_n(x,\vp))+f_n(x)}{n}>s,~\forall~n\geq N,\vp=\frac1{M}\}.
\end{align*}
Thus $Z_{M,N}\subset Z_{M,N+1}$ and $\cup_{N=1}^\infty Z_{M,N}=Z_M$. Again, we can find $N^*\geq1$ such that $\mu(Z_{M,N^*})>\frac1{2}\mu(Z_M)>0$. For every $x\in Z_{M,N^*}$, $n\geq N^*$, and $0<\vp<\frac1{M}$, we have
\[\mu\big(B_n(x,\vp)\big)\leq e^{-ns+f_{n}(x)}.\]
Fix $\vp\in(0,\frac{1}{M})$ and $\delta\in(0,Z_{M,N^*})$. Suppose $\calF=\{B_{n_i}(y_i,\frac\vp2)\}_{i\geq1}$ is a countable family such that $n_i\geq N^*$, $x\in X$ and $\mu\big(\cup_iB_{n_i}(y_i,\frac\vp2)\big)\geq1-\delta$.
Let
\[\calF'=\{B_{n_i}(y_i,\frac\vp2):Z_{M,N^*}\cap B_{n_i}(y_i,\frac\vp2)\neq\emptyset\}.\]
For each $B_{n_i}(y_i,\frac\vp2)\in\calF'$, there exists $x_i\in Z_{M,N^*}\cap B_{n_i}(y_i,\frac\vp2)$. Applying the triangle inequality yields
\[B_{n_i}(y_i,\frac\vp2)\subset B_{n_i}(x_i,\vp).\]
Let $\calF''=\{(n_i,y_i,x_i)\}$.
It follows that
\begin{align*}
\sum_{\calF}e^{-n_is+f_{n_i}(y_i,\frac\vp2)}&\geq\sum_{\calF'}e^{-n_is+f_{n_i}(y_i,\frac\vp2)}\\
&\geq\sum_{\calF''}e^{-n_is+f_{n_i}(x_i)}\\
&\geq\sum_{\calF''}\mu\big(B_{n_i}(x_i,\vp)\big)\\
&\geq\sum_{\calF''}\mu\big(B_{n_i}(y_i,\frac\vp2)\big)\\
&\geq\mu(Z_{M,N^*}\cap\bigcup_{\calF''}B_{n_i}(y_i,\frac\vp2))\\
&=\mu(Z_{M,N^*}\cap\bigcup_{\calF}B_{n_i}(y_i,\frac\vp2))>\mu(Z_{M,N^*})-\delta.
\end{align*}
Thus,
\[\calM_\mu(n,s,\frac\vp2,\delta,T,f)\geq\mu(Z_{M,N^*})-\delta>0.\]
This implies that $\calM_\mu(s,\frac\vp2,\delta,T,f)>0$ and ${P^{KB}_\mu}'(\frac\vp2,\delta,T,f)\geq s$.
Therefore, $P_\mu^{KB}(T,f)\geq s$. Since $s<\underline{P}_u(T,f)$ is arbitrary, the desired inequality follows.


Next, we show the second inequality
\[\overline{P}_u(T,f)\leq P^{KP}_{\mu}(T,f).\]
For any $s<\overline{P}_u(T,f)$, as we did in the proof of the first inequality, we can find $\vp,\beta>0$, and a Borel set $A\subset X$ with $\mu(A)>0$ such that
\[\limsup_{n\to\infty}\frac{-\log\mu(B_n(x,\vp))+f_n(x)}{n}>s+\beta,~\forall x\in A.\]
Fix $\delta\in(0,\mu(A))$. We shall show that $P^{KP}_\mu(\frac{\vp}{5},\delta,T,f)\geq s$, which implies that $P^{KP}_\mu(T,f)\geq s$. It suffices to show that
\[M^{\calP}_\mu(s,\frac{\vp}{5},\delta,T,f)=\infty.\]
Let $\{Z_i\}_{i\geq1}$ be a countable family with $\mu\big(\cup_iZ_i\big)>1-\delta$. Since
\[A=(A\cap \bigcup_{i\in I}Z_i)\cup(X\setminus\bigcup_{i\in I}Z_i),\]
it follows that $\mu\big((A\cap \cup_{i\in I}Z_i)\big)\geq \mu(A)-\delta>0$.
Thus there exists $i$ such that $\mu(A\cap Z_i)>0$. From the proof of the second assertion in Theorem~\ref{thm:Billingsley-theorem}, we have
\[M^{P}(s,\frac{\vp}{5},Z_i,T,f)\geq M^{P}(s,\frac{\vp}{5},A\cap Z_i,T,f)=\infty.\]
This implies that
\[M^{\calP}_\mu(s,\frac{\vp}{5},\delta,T,f)=\infty.\]
\end{proof}

\section{Proofs of Theorems~\ref{thm:variational-p-for-BP} and~\ref{thm:variational-p-for-P}}\label{sec:VP}

\textbf{Proof of Theorem~\ref{thm:variational-p-for-BP}.}
It is clear that
\begin{align*}
\sup\{P_\mu^{B}(T,f):\mu\in{\calM}(X),~\mu(Z)=1\}\leq P^B(Z,T,f).
\end{align*}
Utilizing Proposition~\ref{prop:MKPs-equal-MPs} yields
\[\sup\{P_\mu^{KB}(T,f):\mu\in{\calM}(X),~\mu(Z)=1\}\leq \sup\{P_\mu^{B}(T,f):\mu\in{\calM}(X),~\mu(Z)=1\}.\]
Proposition~\ref{prop:m-p-leq-kp} implies that
\[\sup\{\underline{P}_{\mu}(T,f):\mu\in{\calM}(X),~\mu(Z)=1\}
\leq\sup\{P_\mu^{KB}(T,f):\mu\in{\calM}(X),~\mu(Z)=1\}.\]
Employing Theorem B in~\cite{Tang2015Variational} completes the proof, which tells us that
\[P^B(Z,T,f)=\sup\{\underline{P}_\mu(T,f):\mu\in\calM(X),~\mu(Z)=1\}.\]
\hfill$\Box$

\textbf{Proof of Theorem~\ref{thm:variational-p-for-P}.}
Using Proposition~\ref{prop:m-p-leq-kp} and Proposition~\ref{prop:MKPs-equal-MPs}, this shows that
\begin{align*}
\sup\{\overline{P}_{\mu}(T,f):\mu\in{\calM}(X),~\mu(Z)=1\}
&\leq\sup\{P_{\mu}^{KP}(T,f):\mu\in{\calM}(X),~\mu(Z)=1\}\\
&\leq\sup\{P_{\mu}^{P}(T,f):\mu\in{\calM}(X),~\mu(Z)=1\}\\
&\leq P^P(Z,T,f).
\end{align*}
We shall prove
\begin{equation}\label{eq:packing-p-leq-up-me}
P^P(Z,T,f)\leq\sup\{\overline{P}_{\mu}(T,f):\mu\in{\calM}(X),~\mu(Z)=1\}.
\end{equation}
For this end, we employ the approach used by Feng and Huang in~\cite{feng2012variational}. The following lemma is very important in the proof.

\begin{lem}\label{lem:lemma-for-packing-va}
Let $Z\subset X$, $\vp>0$, and $s>\|f\|_{\infty}$. If $M^{P}(s,\varepsilon,Z,T,f)=\infty$, then for any given finite interval $(a,b)\subset[0,\infty)$ and $N\in\N$, there exists a finite disjoint collection  $\{\overline{B}_{n_i}(x_i,\vp)\}$ such that $x_i\in Z$, $n_i\geq N$, and $\sum_ie^{-n_is+f_{n_i}(x_i)}\in(a,b)$.
\end{lem}

\begin{proof}
Take $N_1>N$ large enough such that  $e^{N_1(\|f\|_\infty-s)}<b-a$. Since $M^{P}(s,\varepsilon,Z,T,f)=\infty$, it follows that $M^{P}(N_1,s,\varepsilon,Z,T,f)=\infty$. There hence exists a finite disjoint collection $\{\overline{B}_{n_i}(x_i,\vp)\}$ such that $x_i\in Z$, $n_i\geq N_1$ and $\sum_ie^{-n_is+f_{n_i}(x_i)}>b$.
Since $e^{-n_is+f_{n_i}(x_i)}\leq e^{n_i(\|f\|_\infty-s)}<b-a$, we can discard elements in this collection one by one until we have $\sum_ie^{-n_is+f_{n_i}(x_i)}\in(a,b)$.
\end{proof}

We now show that inequality~\eqref{eq:packing-p-leq-up-me} holds.
For any $s\in(\|f\|_\infty,P^P(Z,T,f))$. Take $\vp$ small enough such that $s<P^P(\vp,Z,T,f)$. Fix $t\in(s,P^P(\vp,Z,T,f))$. We are going to construct inductively a sequence of finite sets $(K_i)^\infty_{i=1}$ and a sequence of finite measures $(\mu_i)^\infty_{i=1}$ such that $K_i\subset K$ and $\mu_i$ is supported on $K_i$ for each $i$. Together with these two sequences, we construct also a sequence of positive numbers $(\gamma_i)$, and a sequence of integer-valued functions ($m_i: K_i \rightarrow\N$). The construction is divided into three steps:

\emph{Step 1}. Construct $K_1$ and $\mu_1$ as well as $m_1(\cdot)$ and $\gamma_1$.

Note that $M^{\calP}(t,\vp,Z,T,f)=\infty$. Let
\[H=\bigcup\{G\subset X:G~\text{is open},~M^{\calP}(t,\vp,Z\cap G,T,f)=0\}.\]
Then $M^{\calP}(t,\vp,Z\cap H,f)=0$ by the separability of $X$. Let $Z'=Z\setminus H=Z\cap(X\setminus H)$.

\textbf{Claim 1.} For any open set $G\subset X$, either $Z'\cap G=\emptyset$ or $M^{\calP}(t,\vp,Z'\cap G,T,f)>0$.

\textbf{Proof of Claim 1.} Suppose that $M^{\calP}(t,\vp,Z'\cap G,T,f)=0$. Since $Z=Z'\cup(Z\cap H)$,
\[M^{\calP}(t,\vp,Z\cap G,T,f)\leq M^{\calP}(t,\vp,Z'\cap G,T,f)+M^{\calP}(t,\vp,Z\cap H,T,f)=0.\]
Thus $G\subset H$, which implies that $Z'\cap G=\emptyset$.

Since
\[M^{\calP}(t,\vp,Z,T,f)\leq M^{\calP}(t,\vp,Z\cap H,T,f)+M^{\calP}(t,\vp,Z',T,f)\]
and $M^{\calP}(t,\vp,Z\cap H,T,f)=0$,
\[M^{\calP}(t,\vp,Z',T,f)=M^{\calP}(t,\vp,Z,T,f)=\infty.\]
It then follows that
\[M^{\calP}(s,\vp,Z',T,f)=\infty.\]
Using Lemma~\ref{lem:lemma-for-packing-va}, we can find a finite set $K_1\subset Z'$, an integer-valued function $m_1(x)$ on $K_1$ such that the collection $\{\overline{B}_{m_1(x)}(x,\vp)\}_{x\in K_1}$ is disjoint and
\[\sum_{x\in K_1}e^{-m_1(x)s+f_{m_1(x)}(x)}\in(1,2).\]
Define
\[\mu_1=\sum_{x\in K_1}e^{-m_1(x)s+f_{m_1(x)}(x)}\delta_{x},\]
where $\delta_{x}$ denotes the Dirac measure at $x$. Take a small $\gamma_1>0$ such that for any $z\in\overline{B}(x,\gamma_1)$ we have
\begin{equation}\label{eq:pairwise-disjoint-1}
\big(\overline{B}(z,\gamma_1)\cup\overline{B}_{m_1(x)}(z,\vp)\big)\bigcap\big(\bigcup_{y\in K_1\setminus\{x\}}
\overline{B}(y,\gamma_1)\cup\overline{B}_{m_1(y)}(y,\vp)
\big)=\emptyset.
\end{equation}
Since $K_1\subset Z'$, by Claim 1, for any $x\in K_1$,
\[M^{\calP}(t,\vp,Z\cap B(x,\gamma_1/4),T,f)
\geq M^{\calP}(t,\vp,Z'\cap B(x,\gamma_1/4),T,f)>0.\]

\emph{Step 2}. Construct $K_2$ and $\mu_2$ as well as $m_2(\cdot)$ and $\gamma_2$.

By (\ref{eq:pairwise-disjoint-1}), the family of balls $\{\overline{B}(x,\gamma_1)\}_{x\in K_1}$ are pairwise disjoint. For each $x\in K_1$, since $M^{\calP}(t,\vp,Z\cap B(x,\gamma_1/4),f)>0$, we can construct as in Step 1 a finite set
\[E_2(x)\subset Z\cap B(x,\gamma_1/4)\]
and an integer-valued function
\[m_2:E_2(x)\to\N\cap[\max\{m_1(y):y\in K_1\},\infty)\]
such that
\begin{itemize}
  \item[(2-a)] $M^{\calP}(t,\vp,Z\cap G,T,f)>0$ for each open set $G$ with $G\cap E_2(x)\neq\emptyset$;
  \item[(2-b)] The elements in $\{\overline{B}_{m_2(y)}(y,\vp)\}_{y\in E_2(x)}$ are disjoint, and
  \[\mu_1(\{x\})<\sum_{y\in E_2(x)}e^{-m_2(y)s}<(1+2^{-2})\mu_1(\{x\}).\]
\end{itemize}
To see it, we fix $x\in K_1$. Denote $F=Z\cap B(x,\gamma_1/4)$. Let
\[H_x=\bigcup\{G\subset X:G~\text{is open},~M^{\calP}(t,\vp,F\cap G,T,f)=0\}.\]
Set $F'=F\setminus H_x$. Then as in Step 1, we can show that
\[M^{\calP}(t,\vp,F',T,f)=M^{\calP}(t,\vp,F,T,f)>0\]
and
\[M^{\calP}(t,\vp,F'\cap G,T,f)>0\]
for any open set $G$ with $G\cap F'\neq\emptyset$.
Since $s<t$,
\[M^{\calP}(s,\vp,F',T,f)=\infty.\]
Using Lemma~\ref{lem:lemma-for-packing-va} again, we can find a finite set $E_2(x)\subset F'$ and a map $m_2:E_2(x)\to\N\cap[\max\{m_1(y):y\in K_1\},\infty)$ so that (2-b) holds. Observe that if $G\cap E_2(x)\neq\emptyset$ and $G$ is open, then $G\cap F'\neq\emptyset$. Hence
\[M^{\calP}(t,\vp,Z\cap G,T,f)\geq M^{\calP}(t,\vp,F'\cap G,T,f)>0.\]
Thus (2-a) holds.

Since the family $\{\overline{B}(x,\gamma_1)\}_{x\in K_1}$ is disjoint, $E_2(x)\cap E_2(x')=\emptyset$ for different $x,x'\in K_1$. Define $K_2=\cup_{x\in K_1}E_2(x)$ and
\[\mu_2=\sum_{y\in K_2}e^{-m_2(y)s+f_{m_2(y)}(y)}\delta_y.\]
By~\ref{eq:pairwise-disjoint-1} and (2-b), the elements in $\{\overline{B}_{m_2(y)}(y,\vp)\}_{y\in K_2}$ are pairwise disjoint. Hence we can take $\gamma_2\in(0,\gamma_1/4)$ such that for any function $z:K_2\to X$ with $d(x,z(x))\leq\gamma_2$ for $x\in K_2$ , we have
\begin{equation}\label{eq:pairwise-disjoint-2}
\big(\overline{B}(z(x),\gamma_2)\cup\overline{B}_{m_2(x)}(z(x),\vp)\big)\bigcap\big(\bigcup_{y\in K_2\setminus\{x\}}
\overline{B}(z(y),\gamma_2)\cup\overline{B}_{m_2(y)}(z(y),\vp)
\big)=\emptyset
\end{equation}
for each $x\in K_2$.
Since $x\in K_2$, there exists $y\in K_1$ such that $x\in E_2(y)$. By (2-a), for each $x\in K_2$,
\[M^{\calP}(t,\vp,Z\cap B(x,\gamma_2/4),T,f)>0.\]

\emph{Step 3}. Assume that $K_i$, $\mu_i$, $m_i(\cdot)$, and $\gamma_i$ have been constructed for $i=1,\ldots,p$. In particular, suppose that for any function $z:K_p\to X$ wiht $d(x,z(x))<\gamma_p$ for each $x\in K_p$, we have
\begin{equation}\label{eq:pairwise-disjoint-p}
\big(\overline{B}(z(x),\gamma_p)\cup\overline{B}_{m_p(x)}(z(x),\vp)\big)\bigcap\big(\bigcup_{y\in K_p\setminus\{x\}}
\overline{B}(z(y),\gamma_p)\cup\overline{B}_{m_p(y)}(z(y),\vp)
\big)=\emptyset
\end{equation}
for each $x\in K_p$; and $M^{\calP}(t,\vp,Z\cap B(x,\gamma_p/4),f)>0$ for each $x\in K_p$. We construct below $K_{p+1}$, $\mu_{p+1}$, $m_{p+1}(\cdot)$ and $\gamma_{p+1}$ in a way similar to Step 2.

Note that the elements in $\{\overline{B}(x,\gamma_p)\}_{x\in K_p}$ are pairwise disjoint. For each $x\in K_p$, since $M^{\calP}(t,\vp,Z\cap B(x,\gamma_p/4),f)>0$, we can construct as in Step 2, a finite set
\[E_ {p+1}(x)\subset Z\cap B(x,\gamma_p/4)\]
and an integer-valued function
\[m_{p+1}:E_{p+1}(x)\to\N\cap[\max\{m_p(y):y\in K_p\},\infty)\]
such that
\begin{itemize}
  \item[(3-a)] $M^{\calP}(t,\vp,Z\cap G,T,f)>0$ for each open set $G$ wiht $G\cap E_{p+1}(x)\neq\emptyset$;
  \item[(3-b)] The elements in $\{\overline{B}_{m_{p+1}(y)}(y,\vp)\}_{y\in E_{p+1}(x)}$ are disjoint and
  \[\mu_p(\{x\})<\sum_{y\in E_{p+1}(x)}e^{-m_{p+1}(y)s}\leq(1+2^{-p-1})\mu_p(\{x\}).\]
\end{itemize}
Clearly
\[E_{p+1}(x)\cap E_{p+1}(y)=\emptyset\]
for different $x,y\in K_p$.
Define $K_{p+1}=\cup_{x\in K_p}E_{p+1}(x)$ and
\[\mu_{p+1}=\sum_{y\in K_{p+1}}e^{-m_{p+1}(y)s+f_{m_{p+1}}(y)}\delta_y.\]
By~\eqref{eq:pairwise-disjoint-p} and (3-b), the elements in $\{\overline{B}(x,\vp)\}_{x\in K_{p+1}}$ are disjoint. Hence we can take small enough $\gamma_{p+1}\in(0,\gamma_p/4)$ such that for any function $z:K_2\to X$ with $d(x,z(x))\leq\gamma_2$, we have for each $x\in K_{p+1}$,
\begin{equation}\label{eq:pairwise-disjoint-2}
\big(\overline{B}(z(x),\gamma_{p+1})\cup\overline{B}_{m_{p+1}(x)}(z(x),\vp)\big)\bigcap\big(\bigcup_{y\in K_{p+1}\setminus\{x\}}
\overline{B}(z(y),\gamma_{p+1})\cup\overline{B}_{m_{p+1}(y)}(z(y),\vp)
\big)=\emptyset.
\end{equation}
Since $x\in K_{p+1}$, there exists $y\in K_p$ such that $x\in E_{p+1}(y)$. Thus by (2-a), for any $x\in K_{p+1}$,
\[M^{\calP}(t,\vp,Z\cap B(x,\gamma_{p+1}/4),T,f)>0.\]

As in the above steps, we can construct by induction  $\{K_i\}$, $\{\mu_i\}$, $\{m_i(\cdot)\}$, and $\{\gamma_i\}$. We summarize some of their basic properties as follows:

(a) For each $i$, the family $\calF_i:=\{\overline{B}(x,\gamma_i):x\in K_i\}$ is disjoint. For every $B\in\calF_{i+1}$, there exists $x\in K_i$ such that $B\subset \overline{B}(x,\gamma_i/2)$.

(b) For each $x\in K_i$ and $z\in\overline{B}(x,\gamma_i)$,
\begin{equation}\label{eq:1-of-b}
\overline{B}_{m_i(x)}(z,\vp)\cap\bigcup_{y\in K_i\setminus\{x\}}\overline{B}(y,\gamma_i)=\emptyset
\end{equation}
and
\begin{equation}\label{eq:2-of-b}
\mu_i(\overline{B}(x,\gamma_i))=e^{-m_i(x)s}
\leq\sum_{y\in E_{i+1}(x)}e^{-m_{i+1}(y)s+f_{m_{i+1}}(y)}
\leq(1+2^{-i-1})\mu_i(\overline{B}(x,\gamma_i)),
\end{equation}
where $E_{i+1}(x)=\overline{B}(x,\gamma_i)\cap K_{i+1}$.
By~(\ref{eq:2-of-b}), we have
\begin{align*}
\mu_i(F_i)&\leq \mu_{i+1}(F_i)=\sum_{F\in\calF_{i+1}:F\subset F_i}\mu_{i+1}(F)\\
&\leq\sum_{F\in\calF_{i+1}:F\subset F_i}(1+2^{-i-1})\mu_i(F)\\
&=(1+2^{-i-1})\sum_{F\in\calF_{i+1}:F\subset F_i}\mu_i(F)\\
&\leq(1+2^{-i-1})\mu_i(F_i), ~F_i\in\calF_i.
\end{align*}
Using the above inequalities repeatedly, we have for any $j>i$,
\begin{equation}\label{eq:inequality-between-iFi}
\mu_i(F_i)\leq\mu_j(F_i)\leq\prod_{n=i+1}^j(1+2^{-n})\mu_i(F_i)\leq C\mu_i(F_i),~\forall F_i\in\calF_i,
\end{equation}
where $C:=\prod_{n=1}^\infty(1+2^{-n})<\infty$.

Let $\hat{\mu}$ be a limit point of $\{\mu_i\}$ in the weak-star topology. Let
\[K^*=\bigcap_{n=1}^\infty\overline{\bigcup_{i\geq n}K_i}=\lim_{n\to\infty}\overline{\bigcup_{i\geq n}K_i}.\]
Then $\hat{\mu}$ is supported on $K^*$, $K^*\subset Z$, and for every $i\in\N$, $K^*\subset\cup_{x\in K_i}B(x,\gamma_i)$.

By~\eqref{eq:inequality-between-iFi},
\begin{align*}
e^{-m_i(x)s+f_{m_i(x)}(x)}=\mu_i(\overline{B}(x,\gamma_i))
\leq\hat{\mu}(\overline{B}(x,\gamma_i))\leq C\mu_i(\overline{B}(x,\gamma_i))
=Ce^{-m_i(x)s+f_{m_i(x)}(x)},~~\forall x\in K_i.
\end{align*}
In particular,
\[1\leq \sum_{x\in K_1}\mu_1(B(x,\gamma_1))\leq \sum_{x\in K_1}\hat{\mu}(B(x,\gamma_1))=\hat{\mu}(K^*)\leq\sum_{x\in K_1}C\mu_1(B(x,\gamma_1))\leq 2C.\]
By~(\ref{eq:1-of-b}), for every $x\in K_i$ and $z\in\overline{B}(x,\gamma_i)$,
\[\hat{\mu}(\overline{B}_{m_i(x)}(z,\vp))\leq\hat{\mu}(\overline{B}(x,\gamma_i))\leq Ce^{-m_i(x)s+f_{m_i(x)}(x)}.\]
For each $z\in K^*$ and $i\in\N$, there exists $x\in K_i$ such that $z\in\overline{B}(x,\gamma_i)$. Thus
\[\hat{\mu}(\overline{B}_{m_i(x)}(z,\vp))\leq Ce^{-m_i(x)s+f_{m_i(x)}(x)}.\]
Let $\mu=\hat{\mu}/\hat{\mu}(K^*)$. Then $\mu\in\mathcal{M}(X)$, $\mu(K^*)=1$, and for every $z\in K^*$, there exists a sequence $\{k_i\}_{i\geq1}$ with $k_i\to\infty$ such that
\[\mu(B_{k_i}(z,\vp))\leq\frac{ Ce^{-k_is+f_{k_i}(z)}}{\hat{\mu}(K^*)}.\]
This implies that $\overline{P}_{\mu}(T,f)\geq s$.
\hfill$\square$

\section{Proof of Theorem~\ref{thm:va-pp-gene}}\label{sec:va-of-gene}
In this section, we prove Theorem~\ref{thm:va-pp-gene}. Before the proof, we need some lemmas. Recall that for every $x\in X$,
\[\Upsilon_n(x)=\frac{1}{n}\sum_{i=0}^{n-1}\delta_{T^i(x)}.\]
A point $x\in X$ is generic for $\mu$ if $\lim_{n\to\infty}\Upsilon_n(x)=\mu$ in the weak$^\ast$-topology. $G_\mu$ is the set of all generic points for $\mu$.

Let $F\subset\calM(X)$ be a neighborhood of $\mu\in\calM(X,T)$ and $n\in\N$. Define
\[X_{n,F}=\{x\in X:\Upsilon_n(x)\in F\}.\]

Recall that $E\subset X$ is said to be an \emph{$(n,\vp)$-separated subset} if $x, y\in E$, $x\neq y$ implies $d_n(x,y)>\vp$.

Let
\[P(n,\vp,X_{m,F},T,f)=\sup\{\sum_{x\in E}e^{f_n(x)}|E~\text{is an}~(n,\vp)~\text{separated subset of}~X_{m,F}\}.\]

\begin{lem}\label{lem:sepa-leq-mea-f}
Let $(X,T)$ be a TDS, $f\in C(X,\R)$, and ${E_n}$ be a sequence of
$(n,\vp)$-separated subsets. Define
\[\mu_n=\frac{1}{n\sum_{z\in E_n}e^{f_n(z)}}\sum_{y\in E_n}\sum_{i=0}^{n-1}e^{f_n(y)}\delta_{T^i(y)}.\]
Suppose that $\lim_{n\to\infty}\mu_n=\mu$. Then
\[\limsup_{n\to\infty}\frac{1}{n}\log\sum_{y\in E_n}e^{f_n(y)}\leq h_u(T)+\int f\md\mu.\]
\end{lem}

\begin{proof}
The result is contained in the second part of
the proof of Theorem 9.10 in~\cite{Walters1982}.
\end{proof}

\begin{lem}\label{lem:vp-inf-leq-me-f}
Let $(X,T)$ be a TDS and $\mu\in\calM(X,T)$. Then
\[\lim_{\vp\to0}\inf_{\mu\in F}\limsup_{n\to\infty}\frac{1}{n}\log P(n,\vp,X_{n,F},T,f)\leq h_\mu(T)+\int f\md\mu,\]
where the infimum is taken over all neighborhood $F\subset\calM(X)$ of $\mu$.
\end{lem}

\begin{proof}
Here we modify the proof of Proposition 3.1 in~\cite{Pfister-Sullivan2007On}. If $h_\mu+\int f\md\mu=\infty$, then there is nothing to prove. Let $h_\mu+\int f\md\mu<\infty$. Assume that
\[\lim_{\vp\to0}\inf_{\mu\in F}\limsup_{n\to\infty}\frac{1}{n}\log P(n,\vp,X_{n,F},T,f)> h_\mu(T)+\int f\md\mu.\]
Then there exist $\delta>0$ and $\vp^*>0$ such that
\[\inf_{\mu\in F}\limsup_{n\to\infty}\frac{1}{n}\log P(n,\vp,X_{n,F},T,f)\geq h_\mu(T)+\int f\md\mu+2\delta,~\forall \vp<\vp^*.\] Fix $\vp\in(0,\vp^*)$. There exists a decreasing sequence of convex closed neighborhoods $\{C_n\}$ of $\mu$ such that $\cap_nC_n=\{\mu\}$ and
\[\limsup_{n\to\infty}\frac{1}{n}\log P(n,\vp,X_{n,C_n},T,f)\geq h_\mu(T)+\int f\md\mu+2\delta.\]
Let $E_n\subset X_{n,C_n}$ be an $(n,\vp)$ separated subset such that
\[\sum_{x\in E_n}e^{f_n(x)}>\frac{n}{n+1}P(n,\vp,X_{n,C_n},T,f).\]
Then
\begin{align*}
\limsup_{n\to\infty}\frac{1}{n}\log\sum_{x\in E_n}e^{f_n(x)}&\geq\limsup_{n\to\infty}\frac{1}{n}\log\frac{n}{n+1}P(n,\vp,X_{n,C_n},T,f)\\
&=\limsup_{n\to\infty}\frac{1}{n}\log P(n,\vp,X_{n,C_n},T,f)\\
&\geq h_\mu(T)+\int f\md\mu+2\delta.
\end{align*}
Define
\[\mu_n=\frac{1}{n\sum_{z\in E_n}e^{f_n(z)}}\sum_{y\in E_n}\sum_{i=0}^{n-1}e^{f_n(y)}\delta_{T^i(y)}\in C_n.\]
By definition $\lim_{n\to\infty}\mu_n=\mu$. It then follows from Lemma~\ref{lem:sepa-leq-mea-f} that
\[\limsup_{n\to\infty}\frac{1}{n}\log\sum_{y\in E_n}e^{f_n(y)}\leq h_u(T)+\int f\md\mu,\]
which is a contradiction.
\end{proof}

\begin{lem}\label{lem:pp-leq-me-f}
Let $(X,T)$ be a TDS, $\mu\in\calM(X,T)$, $f\in C(X,\R)$, and $G_\mu$ be the set of all generic points for $\mu$. Then
\[P^P(G_\mu,T,f)\leq h_u(T)+\int f\md\mu.\]
\end{lem}

\begin{proof}
For any neighborhood $F\subset\calM(X)$ of $\mu$ and $m\in\N$, let
\[G_\mu^m=\{x\in G_\mu:\Upsilon_n(x)\in F~\forall n\geq m\}.\]
Thus $G_\mu=\cup_{m=1}^\infty G_\mu^m$ and $G_\mu^m\subset X_{n,F}$ for all $n\geq m$. Fix $m\in\N$. Let $\alpha<\beta<P^P(G_\mu^m,T,f)$. Then there exists $\gamma>0$ such that $P^P(\vp,G_\mu^m,T,f)>\beta$ for all $\vp\in(0,\gamma)$. It follows that $M^P(\beta,\vp,G_\mu^m,T,f)\geq M^{\calP}(\beta,\vp,G_\mu^m,T,f)=\infty$. Since $M^P(N,\beta,\vp,G_\mu^m,T,f)$ decreases as $N$ increases, we get $M^P(N,\beta,\vp,G_\mu^m,T,f)=\infty$ for each $N\in\N$. Given $N\geq m$, we can find a finite or countable pairwise disjoint family $\{\overline{B}_{n_i}(x_i,\vp)\}_i$ such that $x_i\in G_\mu^m$, $n_i\geq N$, and
\[\sum_{i}e^{-\beta n_i+f_{n_i}(x_i)}>1,~\forall n_i\geq N.\]
For each $k\geq N\geq m$, let $G_\mu^{m,k}=\{x_i\in G_\mu^m:n_i=k\}$. Then
\[\sum_{k=N}^\infty \big(e^{-\beta k}\sum_{x\in G_\mu^{m,k}}e^{f_k(x)}\big)=\sum_{i}e^{-\beta n_i+f_{n_i}(x_i)}>1.\]
There must be some $k\geq N$ with
\[\sum_{x\in G_\mu^{m,k}}e^{f_k(x)}\geq e^{k\alpha}(1-e^{\alpha-\beta}),\] otherwise the above sum is at most
\[\sum_{k=N}^\infty \big(e^{-\beta k}e^{k\alpha}(1-e^{\alpha-\beta})<1.\]
Since $G_\mu^{m,k}$ is a $(k,\vp)$ separated set of $X_{k,F}$, it follows that
\[P(k,\vp,X_{k,F},T,f)\geq e^{k\alpha}(1-e^{\alpha-\beta}).\]
Thus
\[\limsup_{n\to\infty}\frac{1}{n}\log P(n,\vp,X_{n,F},T,f)\geq\alpha.\]
Since $\vp\in(0,\gamma)$ and $F\ni\mu$ are arbitrary, it follows from Lemma~\ref{lem:vp-inf-leq-me-f} that
\[h_\mu(T)+\int f\md\mu\geq\lim_{\vp\to0}\inf_{\mu\in F}\limsup_{n\to\infty}\frac{1}{n}\log P(n,\vp,X_{n,F},T,f)\geq\alpha.\]
Thus $h_\mu(T)+\int f\md\mu\geq P^P(G_\mu^m,T,f)$. It follows that, by Proposition ~\ref{prop:properties-of-pressures},
\[h_\mu(T)+\int f\md\mu\geq P^P(G_\mu,T,f)=\sup_mP^P(G_\mu^m,T,f).\]
\end{proof}


\textbf{Proof of Theorem~\ref{thm:va-pp-gene}.}
Applying Theorem 3 in~\cite{Pesin1984Topological}, we have $P^B(G_\mu,T,f)=h_u(T)+\int f\md\mu$. The desired equality then follows by employing Lemma~\ref{lem:pp-leq-me-f} and Proposition~\ref{prop:properties-of-pressures}.\hfill$\Box$

\section*{Acknowledgements}
The first author was supported by Guangdong Basic and Applied Basic Research Foundation (No.2019A1515110932);
the second author was supported by the National Natural Science Foundation of China
(Nos.11701584 and 11871228) and the Natural Science Research Project of Guangdong Province (No.2018KTSCX122).


\begin{thebibliography}{10}

\bibitem{Adler1965}
R.~L. Adler, A.~G. Konheim, and M.~H. McAndrew.
\newblock Topological entropy.
\newblock {\em Transactions of the American Mathematical Society},
  114(2):309--319, 1965.

\bibitem{Bowen1971}
R.~Bowen.
\newblock Entropy for group endomorphisms and homogeneous spaces.
\newblock {\em Transactions of the American Mathematical Society},
  153:401--414, 1971.

\bibitem{Bowen1973Topological}
R.~Bowen.
\newblock Topological entropy for noncompact sets.
\newblock {\em Transactions of the American Mathematical Society},
  184:125--136, 1973.

\bibitem{Dinaburg1970}
E.~I. Dinaburg.
\newblock A correlation between topological entropy and metric entropy.
\newblock {\em Dokl.akad.nauk Sssr}, 190:19--22, 1970.

\bibitem{Dou2017Topological}
D.~Dou, M.~Fan, and H.~Qiu.
\newblock Topological entropy on subsets for fixed-point free flows.
\newblock {\em Discrete \verb'&' Continuous Dynamical Systems},
  37(12):6319--6331, 2017.

\bibitem{Dou2020Packing}
D.~Dou, D.~Zheng, and X.~Zhou.
\newblock Packing topological entropy for amenable group actions.
\newblock {\em Ergodic Theory and Dynamical Systems}, pages 1--35, 2021.

\bibitem{feng2012variational}
D.-J. Feng and W.~Huang.
\newblock Variational principles for topological entropies of subsets.
\newblock {\em Journal of Functional Analysis}, 263(8):2228--2254, 2012.

\bibitem{Huang2020A}
X.~Huang, Z.~Li, and Y.~Zhou.
\newblock A variational principle of topological pressure on subsets for
  amenable group actions.
\newblock {\em Discrete \verb'&' Continuous Dynamical Systems},
  40(5):2687--2703, 2020.

\bibitem{Katok1980Lyapunov}
A.~Katok.
\newblock Lyapunov exponents, entropy and periodic orbits for diffeomorphisms.
\newblock {\em Publications Math{\'e}matiques de l'Institut des Hautes
  {\'E}tudes Scientifiques}, 51(1):137--173, 1980.

\bibitem{Kolmogorov1958A}
A.~N. Kolmogorov.
\newblock A new metric invariant of transient dynamical systems and
  automorphisms in lebesgue spaces.
\newblock {\em Dokl. Akad. Nauk SSSR}, 951(5):861--864, 1958.

\bibitem{Kong2014Slow}
D.~Kong and E.~Chen.
\newblock Slow entropy for noncompact sets and variational principle.
\newblock {\em Journal of Dynamics \verb'&' Differential Equations},
  26(3):477--492, 2014.

\bibitem{Liang2021Packing}
R.~Liang and H.~Lei.
\newblock Packing entropy for fixed-point free flows.
\newblock {\em arXiv:2102.10281}, 2021.

\bibitem{Ma-Wen2008A-Billingsley}
J.-H. Ma and Z.-Y. Wen.
\newblock A \uppercase{B}illingsley type theorem for \uppercase{B}owen entropy.
\newblock {\em Comptes Rendus Mathematique}, 346(9):503--507, 2008.

\bibitem{Mattila1995Geometry}
P.~Mattila.
\newblock {\em Geometry of sets and measures in Euclidean spaces}.
\newblock Cambridge University Press, 1995.

\bibitem{Pesin1997Dimension}
Y.~B. Pesin.
\newblock {\em \uppercase{D}imension \uppercase{T}heory in
  \uppercase{D}ynamical \uppercase{S}ystems: \uppercase{C}ontemporary
  \uppercase{V}iews and \uppercase{A}pplications}.
\newblock University of Chicago Press, Chicago, 1997.

\bibitem{Pesin1984Topological}
Y.~B. Pesin and B.~S. Pitskel.
\newblock Topological pressure and the variational principle for noncompact
  sets.
\newblock {\em Functional Analysis and Its Applications}, 18(4):307--318, 1984.

\bibitem{Pfister-Sullivan2007On}
C.-E. Pfister and W.~G. Sullivan.
\newblock On the topological entropy of saturated sets.
\newblock {\em Ergodic Theory and Dynamical Systems}, 27(3):929--956, 2007.

\bibitem{Ruelle1973statistical}
D.~Ruelle.
\newblock Statistical mechanics on a compact set with $z^\nu$ action satisfying
  expansiveness and specification.
\newblock {\em Transactions of the American Mathematical Society},
  185:237--251, 1973.

\bibitem{Tang2015Variational}
X.~Tang, W.~C. Cheng, and Y.~Zhao.
\newblock Variational principle for topological pressures on subsets.
\newblock {\em Journal of Mathematical Analysis \verb'&' Applications},
  424(2):1272--1285, 2015.

\bibitem{Walters1982}
P.~Walters.
\newblock {\em An introduction to ergodic theory. [Graduate texts in
  mathematics, Vol. 79]}.
\newblock Springer-Verlag, New York, 1982.

\bibitem{Wang2012Variational}
C.~Wang and E.~Chen.
\newblock Variational principles for \uppercase{BS} dimension of subsets.
\newblock {\em Dynamical Systems}, 27(3):359--385, 2012.

\bibitem{Wang2021some}
T.~Wang.
\newblock Some notes on topological and measure-theoretic entropy.
\newblock {\em Qualitative Theory of Dynamical Systems}, 20(1):1--13, 2021.

\bibitem{Wang2020Bowen}
Y.~Wang, E.~Chen, Z.~Lin, and T.~Wu.
\newblock Bowen entropy of sets of generic points for fixed-point free flows.
\newblock {\em Journal of Differential Equations}, 269(11):9846--9867, 2020.

\bibitem{Xu2018Variational}
L.~Xu and X.~Zhou.
\newblock Variational principles for entropies of nonautonomous dynamical
  systems.
\newblock {\em Journal of Dynamics \verb'&' Differential Equations},
  30(3):1053--1062, 2018.

\bibitem{Zheng2016Bowen}
D.~Zheng and E.~Chen.
\newblock Bowen entropy for actions of amenable groups.
\newblock {\em Israel Journal of Mathematics}, 212:895--911, 2016.

\bibitem{Zheng2018Topological}
D.~Zheng and E.~Chen.
\newblock Topological entropy of sets of generic points for actions of amenable
  groups.
\newblock {\em Science China Mathematics}, 61(5):869--880, 2018.

\bibitem{Zhong2020Variationalfree}
X.~F. Zhong and Z.~J. Chen.
\newblock Variational principle for topological pressure on subsets of free
  semigroup actions.
\newblock {\em Acta Mathematica Sinica, English Series}, 37(9):1401--1414,
  2021.

\end{thebibliography}
\end{document}